\documentclass{siamltex}

\usepackage{graphicx}
\usepackage{subfigure}
\usepackage{multirow}
\usepackage[T1]{fontenc}
\usepackage[sc]{mathpazo}
\usepackage{xspace}
\usepackage{url}

\usepackage{amsmath}
\usepackage{amssymb,amsfonts}
\usepackage{upgreek}
\usepackage[hmargin=1.5in]{geometry}

\usepackage{hyperref}
\usepackage{xcolor}
\usepackage{color}
\usepackage{todonotes}
\usepackage{textcomp}
\newcommand{\gray}{\textcolor{gray}}

\graphicspath{{./figures/}}

\usepackage{listings}
\usepackage{caption}
\DeclareCaptionFont{white}{\color{white}}
\DeclareCaptionFormat{listing}{%
  \parbox{\textwidth}{\colorbox{gray}{\parbox{\textwidth}{#1#2#3}}\vskip-4pt}}
\begin{document}
\def\petsc{PETSc\xspace}
\def\numpy{{\tt\bf numpy}\xspace}
\def\clawpack{Clawpack\xspace}
\def\sharpclaw{SharpClaw\xspace}
\def\petscforpy{{\tt\bf petsc4py}\xspace}
\def\pyclaw{{\tt\bf PyClaw}\xspace}
\def\petclaw{{\tt\bf PetClaw}\xspace}
\def\class#1{{\tt\bf #1}\xspace}
\def\function#1{{\tt #1}\xspace}

\newcommand{\Oop}{{\cal O}}
\newcommand{\Dt}{\Delta t}
\newcommand{\Dx}{\Delta x}
\newcommand{\Dy}{\Delta y}
\newcommand{\imh}{{i-\frac{1}{2}}}
\newcommand{\iph}{{i+\frac{1}{2}}}
\newcommand{\bq}{\mathbf q}
\newcommand{\bx}{\mathbf x}
\newcommand{\bflux}{\mathbf f}
\newcommand{\bs}{\mathbf s}

\belowcaptionskip=-10pt

\title{PyClaw: Accessible, Extensible, Scalable Tools for Wave Propagation Problems}
\author{
  David I. Ketcheson\thanks{King Abdullah University of Science and Technology,
    Box 4700, Thuwal, Saudi Arabia, 23955-6900 (\mbox{david.ketcheson@kaust.edu.sa}) }
    \and
   Kyle T. Mandli\thanks{University of Texas at Austin, 1 University Station C0200 Austin, TX 78712-0027 (\mbox{kyle@ices.utexas.edu})}
   \and
   Aron J. Ahmadia\thanks{King Abdullah University of Science and Technology, Box 4700, Thuwal, Saudi Arabia, 23955-6900 (\mbox{aron.ahmadia@kaust.edu.sa}) }
   \and
   Amal Alghamdi\thanks{King Abdullah University of Science and Technology, Box 4700, Thuwal, Saudi Arabia, 23955-6900 (\mbox{amal.alghamdi@kaust.edu.sa}) }
   \and
   Manuel Quezada de Luna\thanks{Texas A\&M University, College Station, TX 77843  (\mbox{mquezada@math.tamu.edu})}
   \and
   Matteo Parsani\thanks{King Abdullah University of Science and Technology, Box 4700, Thuwal, Saudi Arabia, 23955-6900 (\mbox{matteo.parsani@kaust.edu.sa}) }
   \and
   Matthew G. Knepley\thanks{University of Chicago, 5735 S. Ellis Ave. Chicago, IL 60637 (\mbox{knepley@ci.uchicago.edu}) }
   \and Matthew Emmett\thanks{University of North Carolina at Chapel Hill, Chapel Hill, NC 27599 (\mbox{memmett@unc.edu})}}

\maketitle

\begin{abstract}
Development of scientific software involves tradeoffs between ease of use,
generality, and performance.  We describe the design of a general hyperbolic PDE
solver that
can be operated with the convenience of MATLAB yet achieves efficiency near
that of hand-coded Fortran and scales to the largest supercomputers.
This is achieved by using Python for most of the code while employing
automatically-wrapped Fortran kernels for computationally intensive routines,
and using Python bindings to interface with a
parallel computing library and other numerical packages.
The software described here is PyClaw,
a Python-based structured grid solver for general
systems of hyperbolic PDEs \cite{pyclaw}. PyClaw provides a powerful and intuitive interface to the
algorithms of the existing Fortran codes \clawpack and \sharpclaw,
simplifying code development and use while providing
massive parallelism and scalable solvers via the PETSc library.
The package is further
augmented by use of PyWENO for generation of efficient high-order weighted
essentially non-oscillatory reconstruction code.
The
simplicity, capability, and performance of this approach are demonstrated through application to
example problems in shallow water flow, compressible flow and elasticity.
\end{abstract}

\section{Introduction}

Traditionally, scientific codes have been developed in compiled languages
like Fortran or C. There exists an abundance of well-tested,
efficient, often serial implementations of numerical algorithms
in those languages.
It is often desirable to parallelize and extended such codes with new algorithmic
components in order to apply them to increasingly challenging problems.

More recently, high-level scientific programming languages,
such as MATLAB, IDL, and R,
have also become an important platform for numerical codes.  These languages
offer powerful advantages:
they allow code to be written in a language more familiar to scientists
and they permit development to occur in an evolutionary fashion.
Problem parameters can be specified and plotting can be performed
interactively, bypassing the comparatively slow edit/compile/run/plot cycle of
development in Fortran or C \cite{perez2011}.
However, programs written in such high-level languages are not portable to
high performance computing platforms and may be very slow compared to
compiled code.

We present one approach to leveraging the advantages of both kinds of languages.
Our starting point is \clawpack \cite{clawpack45}: a widely used, state-of-the-art
package for solving hyperbolic systems of partial differential equations,
such as those arising in fluid mechanics, astrophysics, geodynamics, magnetohydrodynamics,
oceanography, porous media flow, and numerous other fields.  In this
paper we present PyClaw, a package that brings greater accessibility,
flexibility, and parallelization to Clawpack and related packages.
PyClaw is used as an illustration
to describe, demonstrate, and provide support for a particular approach to construction
of scientific software.  The approach we advocate consists of three steps:
\begin{enumerate}
\item use Python to create a convenient interface to serial legacy code;
\item use Python to interface the resulting code with software tools for
        parallelization of the code, with minimal modification of the serial legacy code;
\item use Python to interface the resulting parallel code to other packages
        that provide extended functionality.
\end{enumerate}
In the case of PyClaw, 
(i) consists of a Python interface to the Clawpack and SharpClaw Fortran-based
packages for numerical solution of systems of hyperbolic PDEs.
This interface allows the code to be operated in
the same convenient and interactive way that one works with MATLAB.
In step (ii), PyClaw was parallelized by interfacing with PETSc, a state-of-the-art 
library for parallel scientific computing.  
This enables parallel computation that scales to large supercomputers and
achieves on-core performance close to that of the legacy code.
Finally, step (iii) is illustrated
in that PyClaw was interfaced with PyWENO, to increase the available order of accuracy of
numerical approximation.  A key in all three steps is the use of the numerical
Python package numpy \cite{numpy}.
The idea of using a layer of numpy-based Python code on top of Fortran, C,
or C++ kernels to solve PDEs has become increasingly popular over the
past several years; see for instance \cite{mardal2007using,FiPy:2009}.
We consider PyClaw to be a very convincing case study.

PyClaw is one of the most highly scalable Python-based codes available.
Perhaps the first well-known scientific project to provide a parallel solver in Python is GPAW,
which extends the Python interpreter itself with parallel data structures and algorithms for electronic structure
calculations \cite{Mortensen2005}.  Python has also previously been used as a tool for parallelizing Fortran and C
codes by introducing parallel code in a Python layer that also calls the
Fortran/C kernels \cite{Nilsen2010}.  In the FiPy package, parallelism is
achieved by using an existing parallel library (Trilinos)
through its Python bindings \cite{FiPy:2009}.  PyClaw takes an approach similar to that of FiPy, in which
all parallel operations over distributed-memory processes are handled by the PETSc library through the
petsc4py Python package
(\url{http://code.google.com/p/petsc4py/}).
This approach offers the advantages of utilizing a documented, well-designed
abstract parallel interface for developers that is already known to
achieve excellent scaling on many architectures.

The algorithms of \clawpack and its high-order extension,
\sharpclaw, are described in Section \ref{sec:hyp}.
The PyClaw framework is described in Section \ref{sec:pyclaw}.
Section \ref{sec:petclaw} describes the parallelization of PyClaw.
As demonstrated in Section
\ref{sec:performance}, PyClaw maintains both the serial performance
of \clawpack and the parallel scalability of PETSc.
In Section \ref{sec:software}, we briefly highlight some of the software
development practices that have contributed to the success of PyClaw.
The combination of wave propagation algorithms and scalable parallelization
enables efficient solution of interesting scientific problems, as demonstrated
through three examples in Section \ref{sec:apps}.  

A repository containing the data, software, and hardware environments for reproducing all
experimental results and figures in this paper is available online (\url{http://bitbucket.org/ahmadia/pyclaw-sisc-rr}).
The most recent release of the PyClaw code is hosted at \url{http://github.com/clawpack/pyclaw} and can be installed
alongside its dependencies in the clawpack distribution in a few minutes with the following \texttt{pip} commands (pip is a freely
available Python package installer and manager):

\begin{lstlisting}[caption=Installing the most recent release of PyClaw and its dependencies]
pip install numpy
pip install clawpack
\end{lstlisting}

The petsc4py package is not an explicit dependency of PyClaw, but if it has been installed, the PetClaw parallel
extension described in Section \ref{sec:petclaw} is seamlessly enabled.

It is our hope that readers will download and try the code.  To whet the reader's appetite,
an example of a complete PyClaw program is shown in Listing \ref{hello-world}.  This
example sets up, runs, and plots the solution of a two-dimensional inviscid fluid dynamics problem
(specifically, test case 6 of \cite{Lax1998}).

\begin{figure}
\begin{lstlisting}[language=python,label=hello-world,caption=Solution of a 2D Euler Riemann problem]
from clawpack import pyclaw
from clawpack import riemann
# Define the equations and the solution algorithm
solver = pyclaw.ClawSolver2D(riemann.rp2_euler_4wave)
# Define the problem domain
domain = pyclaw.Domain([0.,0.],[1.,1.],[100,100])
solver.all_bcs = pyclaw.BC.extrap
solution = pyclaw.Solution(solver.num_eqn,domain)
# Set physical parameters
gamma = 1.4
solution.problem_data['gamma']  = gamma
# Set initial data
xx,yy = domain.grid.p_centers
l = xx<0.5; r = xx>=0.5; b = yy<0.5; t = yy>=0.5
solution.q[0,...] = 2.*l*t + 1.*l*b + 1.*r*t + 3.*r*b
solution.q[1,...] = 0.75*t - 0.75*b
solution.q[2,...] = 0.5*l  - 0.5*r
solution.q[3,...] = 0.5*solution.q[0,...]* \
                    (solution.q[1,...]**2+solution.q[2,...]**2) \
                    + 1./(gamma-1.)
# Solve
solver.evolve_to_time(solution,tend=0.3)            
# Plot
pyclaw.plot.interactive_plot()
\end{lstlisting}
\end{figure}


\section{Finite Volume Hyperbolic PDE solvers\label{sec:hyp}}

    The numerical methods in PyClaw compute approximate solutions of
    systems of hyperbolic conservation laws:
    \begin{align} \label{eq:conslaw}
      \kappa(\bx) \bq_t + \nabla\cdot\bflux(\bq,\bx)_x & = \bs(\bq,\bx).
    \end{align}
    Here $\bq(\bx,t)$ is a vector of conserved quantities (e.g., density,
    momentum, energy) and $\bflux(\bq,\bx)$ represents the flux (modeling wave-like phenomena),
    while $\bs(\bq,\bx)$ represents additional non-hyperbolic {\em source} terms, such as diffusion
    or chemical reactions.  The {\em capacity function} $\kappa(\bx)$ is frequently
    useful for taking into account variations in material properties or
    in the use of non-uniform grids (see \cite[Chapter 6]{levequefvmbook}).
    Here we describe high-resolution shock capturing methods, which
    is one of the most successful classes of numerical methods for solving \eqref{eq:conslaw}.

    Computing solutions to nonlinear hyperbolic equations is often costly.
    Solutions of \eqref{eq:conslaw} generically develop
    singularities (shocks) in finite time, even if the initial data are
    smooth. Accurate modeling of solutions with shocks or strong convective character requires computationally
    expensive techniques, such as Riemann solvers and nonlinear limiters.

    In a finite volume method, the unknowns at time level $t^n$ are taken to be the averages of $q$
    over each cell:
    \begin{align}
        Q^n_i = \frac{1}{\Dx} \int_{x_\imh}^{x_\iph} q(x,t^n) \, dx,
    \end{align}
    where $\Dx = x_\iph - x_\imh$ and $i$ are the local grid spacing and the
    cell's index, respectively.
    A simple update of the cell averages based on the resulting waves gives the classic
    Godunov method, a robust but only first-order accurate numerical scheme:
    \begin{align} \label{eq:fluxdiff}
        Q^{n+1}_i & = Q^n_i - \frac{\Dt}{\Dx}\left(F^n_\iph-F^n_\imh\right),
    \end{align}
    where $F$ and $\Dt$ are the {\em numerical flux} function and the time step. Godunov's
    method results from taking a particular choice of $F$ referred
    to as the {\em upwind flux}.

    The first-order method just described is very dissipative.  Higher-order
    extensions require the computation of higher derivatives of the solution
    or flux.  Near a solution
    discontinuity, or shock, spurious oscillations tend to arise
    due to dispersion and numerical differencing across the discontinuity.
    In order to combat this, shock-capturing methods use special nonlinear
    algorithms to compute numerical derivatives in a non-oscillatory way by
    limiting the the value of the computed derivative in the vicinity of
    a discontinuity \cite{levequefvmbook}. 


    The classic Clawpack algorithm is based on the second-order Lax-Wendroff difference scheme
    that was later extended by LeVeque \cite{leveque1997,levequefvmbook}.
    This scheme can be written in the
    flux-differencing form \eqref{eq:fluxdiff}
    by an appropriate choice of numerical flux, which has the form
    \begin{align} \label{eq:lwflux}
        F^n_\imh & = F_\textup{upwind} + F_\textup{correction},
    \end{align}
    where $F_\textup{upwind}$ is the Godunov flux.
    The classic Clawpack algorithm is based on modifying \eqref{eq:lwflux} by applying
    a limiter to $F_\textup{correction}$. However, in LeVeque's extension, rather than applying the limiter
    to the flux variables, the limiter is applied directly to the
    waves computed by the Riemann solver.  This allows for better accuracy
    by limiting only the characteristic families that are discontinuous
    in a given neighborhood.  Furthermore, the first-order contribution is
    written in terms of {\em fluctuations} (which approximate the quasilinear term $Aq_x$) rather
    than fluxes (which approximate $f(q)$).  This allows the algorithm to be
    applied to hyperbolic systems not in conservation form.

    While the Lax-Wendroff approach can be extended to even higher order,
    this is cumbersome because of the large number of high order terms
    appearing in the Taylor series.  A simpler alternative is the method of
    lines, in which the spatial derivatives are discretized first, leading
    to a system of ODEs that can be solved by traditional methods.  This
    is the approach taken in \sharpclaw \cite{sharpclaw,ketcheson2006,Ketcheson2011}.
    First, a non-oscillatory approximation of the solution is reconstructed
    from the cell averages to give high order accurate point values
    just to the left and right
    of each cell interface.  This reconstruction is performed
    using weighted essentially non-oscillatory (WENO) reconstruction in order to avoid
    spurious oscillations near discontinuities.  
    As in the classic Clawpack algorithm, the scheme is written in terms of fluctuations
    rather than fluxes, so that it can be applied to non-conservative problems.

    \subsection{\clawpack and SharpClaw\label{sub:clawpack}}
    We now describe the "legacy" Fortran codes on which PyClaw is built.

    The classic algorithm implemented in \clawpack (``Conservation Laws Package'') includes
    extensions to two and three dimensions, adaptive mesh refinement, and other
    enhancements \cite{clawpack45}.
    The unsplit multidimensional \clawpack algorithms include additional
    correction terms, computed by a secondary ``transverse'' Riemann solver,
    which approximates corner transport terms.
    The \clawpack software (\url{http://www.clawpack.org}) and its extensions,
    consisting of open source Fortran code, have been freely available since
    1994. More than 7,000 users have registered to download \clawpack.

    \clawpack is a very general tool in the sense that it is easily adapted
    to solve any hyperbolic system of conservation laws.
    The only specialized code required in order to solve a
    particular hyperbolic system is the Riemann solver routine.  A wide range of
    Riemann solvers, including several for the most widely studied hyperbolic systems, have
    been developed by \clawpack users and are also freely available.
    \clawpack handles not only simple Cartesian grids but any logically
    quadrilateral grid, provided that there is a map from a
    uniform Cartesian domain to the desired physical domain
    (see Section \ref{sub:shallow_sphere}).
    Non-hyperbolic source terms ($\bs(\bq,\bx)$) can be easily included via operator splitting.
    For more examples and details regarding Clawpack,
    see \cite{leveque1997} and \cite[Chapter 23]{levequefvmbook}.


    The high-order WENO-based wave propagation method 
    is implemented in \sharpclaw, another Fortran package designed
    similarly to \clawpack and which makes use of the same Riemann solver routines.
    The default options in \sharpclaw employ fifth-order WENO reconstruction
    in space and the fourth-order strong stability preserving (SSP) Runge--Kutta
    method of \cite{ketcheson2008} in time.
    In multi-dimensions \sharpclaw requires propagation of waves only in the
    normal direction to each edge.

\section{PyClaw\label{sec:pyclaw}} 
    PyClaw is an object-oriented framework that incorporates the functionality 
    of \clawpack and \sharpclaw. This functionality may be provided via either 
    pure Python or calls made to the underlying Fortran routines included in 
    \clawpack and \sharpclaw. PyClaw is designed such that the user provides 
    appropriate call-back functions, such as a Riemann solver, that describe 
    and implement the problem in question, PyClaw than manages the ``main'' 
    routine including appropriate time stepping and output. It also avoids the 
    need to deal with strictly formatted data files and reduces the need to 
    write custom Fortran routines for new problems. Instead, problems can be 
    set up interactively or in a simple scripting language. PyClaw also allows 
    for simulation and visualization to be done in a single, interactive 
    environment. Users may engage with this framework at different levels, 
    depending on their expertise and requirements. These interfaces are 
    described in Section \ref{sec:interfaces}.



    PyClaw wraps the full functionality of the ``classic'' 1D, 2D and 3D \clawpack code
    including the use of capacity functions, mapped
    grids, and both dimensionally-split and fully-multidimensional algorithms.
    It also provides the full functionality of SharpClaw, and adds to this with
    higher order WENO reconstruction.
    It does not presently include adaptive mesh refinement, which is part of
    the AMRClaw and GeoClaw extensions of \clawpack.


    \subsection{Interfaces\label{sec:interfaces}}
    The PyClaw distribution includes pre-written application scripts that solve problems
    in acoustics, elasticity, compressible flow, shallow water flow, and other application domains.
    These application scripts represent the typical ``main'' routines that
    lower-level language developers are used to, and are written with ease of
    understanding as their most important goal. 
    These scripts run both on serial workstations and from batch
    processing queues for, e.g., 8,000-node parallel jobs without modification.
    Novice users can solve a wide range of problems by modifying the
    example scripts to deal with different domains, initial conditions,
    boundary conditions, and so forth.  This requires only simple understanding
    of high-level scripting code, but allows users to compute solutions of complex problems
    on large supercomputers.  The scripts for the applications in this paper have
    all been added to the distributed examples, and we plan to continue this practice
    with respect to future publications that use PyClaw.

    Novice or advanced users may also run problems and analyze results in an interactive shell.  When
    Python is invoked with standard input connected to a \verb!tty! device, it reads and executes commands interactively.  This
    feature simplifies serial development, debugging, and visualization, and is familiar to users of commercial software
    such as MATLAB and Mathematica. PyClaw's top level classes present the same API whether used in serial or parallel,
    allowing users to develop interactively, then run production jobs in batch environments.

    Advanced users may want to solve a hyperbolic system that is not included
    among the example applications.  In \clawpack a number of Riemann solvers
    have been written for specific systems and are included with \clawpack.
    Since PyClaw and \clawpack Riemann solvers are interoperable, many systems
    of interest have already been implemented and can be used immediately in
    PyClaw.  A user also has the option of writing his/her own Riemann solver in
    Fortran which can then be utilized in \clawpack as well as PyClaw.
    \clawpack users who wish to run an existing serial \clawpack application in
    parallel can do so easily by wrapping any problem-specific Fortran routines
    (such as those that set initial conditions or compute source terms)
    automatically with f2py and using the resulting Python function handles in
    PyClaw.

    Numerical analysts are often interested in comparing solutions of a problem
    obtained with different methods or different values of method parameters.
    The sample application scripts provide a common functional interface that can be accessed
    from the command line for selecting between solvers, choosing
    serial or parallel computation, and other options.  Users are free to extend this interface to allow more programmatic flexibility at the
    command line or from within a batch environment.  PyClaw also fully exposes the full range of command line options
    available from the underlying PETSc library, allowing advanced users to tweak low-level settings such as message-passing communication strategies.

    Frequently, research scientists are
    interested in comparing the performance of numerical methods.  PyClaw
    enables this comparison by allowing scientific developers to extend the
    software with a new Solver class.  The Solver class, described in the next
    section, is responsible for prescribing a single time step of the numerical
    algorithm.  In PyClaw this is accomplished by implementing a
    \verb!homogenous_step! routine which evolves the solution of the PDE from
    time $t$ to $t+\Delta t$.  Since
    most existing codes have such a routine already, it is often straightforward
    to include legacy code in the PyClaw framework by simply wrapping this
    function. Non-hyperbolic terms $\bs(\bq,\bx)$ can also be incorporated via
    operator splitting.
    This is the most common way to extend the Solver class and
    allows for easy comparison between different numerical methods.

    \subsection{Classes\label{sec:classes}}
    The primary abstractions used in PyClaw are the Solver and the Solution.  The Solver class provides an
    abstract interface to evolve a Solution object forward in time.  The Solution class is a data abstraction
    containing information about the domain of the PDE and the state of the solution at a particular time inside
    the domain.  Here we will discuss these classes and how they interact.

\begin{figure}[ht] 
        \centering
        \subfigure[Structure of a PyClaw Solution object, which may contain multiple Grid objects, each
        of which may have multiple associated State objects.  Each State object has a associated
        fields of conserved quantities (e.g., density, momentum) and optionally, associated
        auxiliary property fields.\label{fig:solution}]{\includegraphics[width=0.45\textwidth]{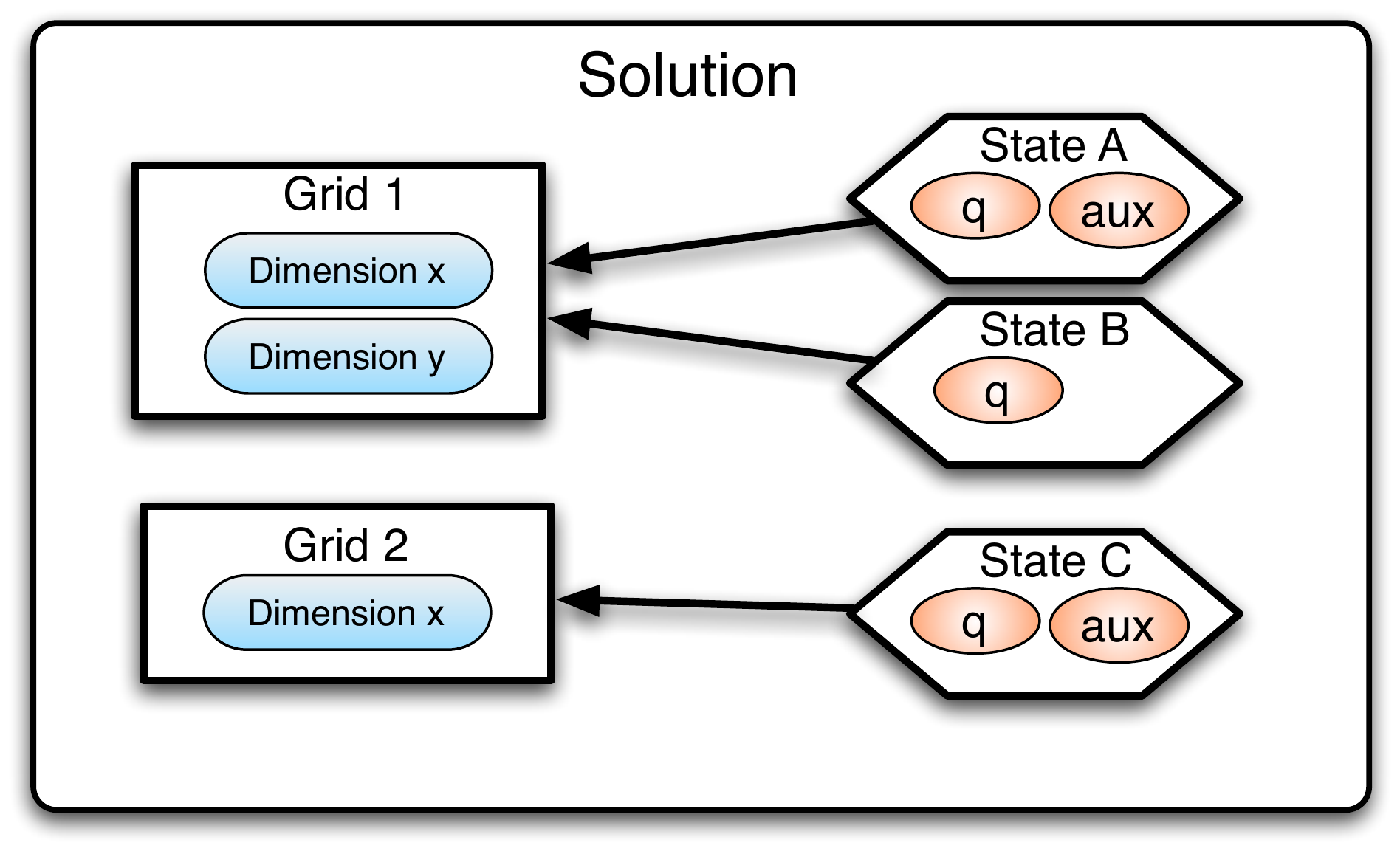}}~~~
        \subfigure[Role of Solver object.  The Solver acts on a Solution object in order to
        evolve it to a later time.\label{fig:solver}]{\includegraphics[width=0.45\textwidth]{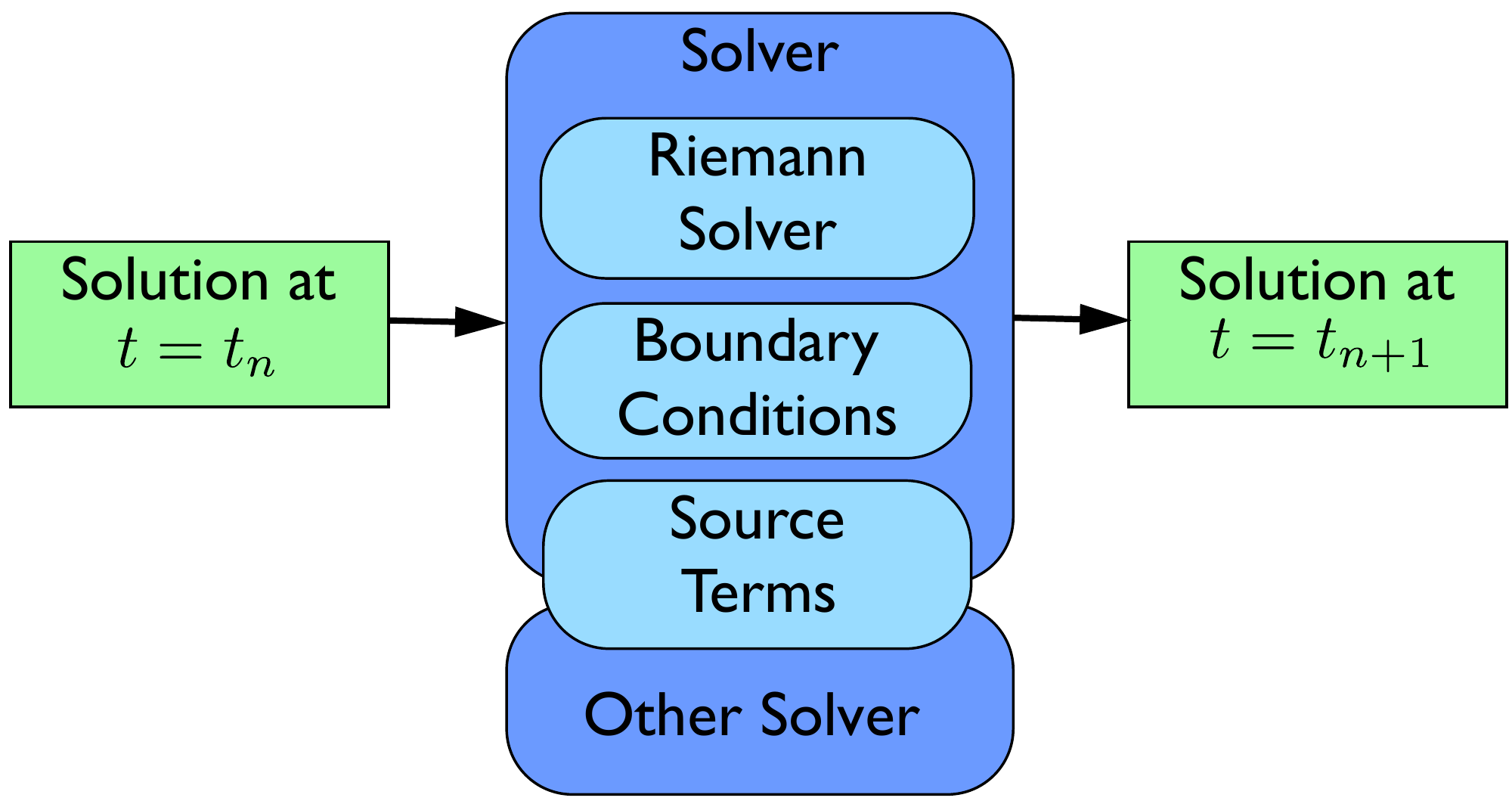}}
        \caption{Structure and function of the main PyClaw classes\label{fig:pyclaw_structure}.}
    \end{figure}

    The role of the Solver is illustrated in Figure~\ref{fig:solver}.
    The Solver class prescribes how a State $Q^n$ is evolved forward in time to obtain $Q^{n+1}$;
    in general this consists of three parts:
    \begin{enumerate}
      \item set ghost cell values based on the prescribed boundary conditions;
      \item advance the solution based on the hyperbolic terms (i.e., $\bq_t + \nabla\cdot\bflux(\bq,\bx) = 0$);
      \item advance the solution based on the source term $\bs(\bq,\bx)$ (if present)
            by a fractional-step approach \cite[Chapter 17]{levequefvmbook}.
    \end{enumerate}
    The base Solver class implements the basic interface to each of these
    functions and a subclass of the Solver class is expected to implement the
    appropriate functions depending on the numerical
    method being implemented.
    The Solver class is sufficiently abstract to accommodate algorithms
    that are based on the method of lines, such as SharpClaw, as well as
    algorithms that are not, such as \clawpack.

    The Solution class, depicted in Figure \ref{fig:solution}, has two purposes:
    \begin{itemize}
        \item Describe the problem domain
        \item Keep track of the values of the state variables $Q^n$ and PDE coefficients
    \end{itemize}
    Like Clawpack, PyClaw simulations are always based on a computational domain
    composed of tensor products of one-dimensional equispaced discretizations of space.  More
    general physical domains may be used as long as it is possible to map them to that computational
    domain.  PyClaw includes a set of geometry classes that implement these abstractions.

    The solution values $Q^n$ and the PDE coefficients are contained in
    numpy arrays stored in the Solution object.
    Thus the full Solution class represents a snapshot of the gridded data. The class acts as a container object with
    one or more Grid and State objects such as in the case of adaptive mesh refinement or nested
    grids, both of which are possible with \clawpack algorithms, though not yet
    available in PyClaw. This hierarchical class
    structure allows the underlying data and algorithms to be modified without the knowledge of the
    interacting objects and without changing the interface presented to the user.  An example of this is the
    PetClaw State object, which reimplements State functionality over distributed memory.  In the future, we intend to
    provide State objects that implement other strategies for accelerating performance.

    \subsection{Extension using PyWENO\label{sub:pyweno}}
    One of the principal achievements of PyClaw has been to facilitate
    the extension of Clawpack and SharpClaw by interfacing with other
    packages.  For example, PyWENO \cite{pyweno} has been used to add
    much higher-order functionality to the existing \sharpclaw code,
    within PyClaw.

    The Fortran code \sharpclaw contains only fifth-order WENO routines
    for evaluation at cell interfaces.
    New WENO routines for PyClaw were generated by PyWENO, which is
    a standalone package for building custom WENO codes.  For a given
    (arbitrarily high) order of reconstruction, PyWENO symbolically
    computes the smoothness coefficients, reconstruction coefficients,
    and optimal (linear) weights of the WENO method.  From these
    symbolic results, PyWENO generates Fortran kernels that perform
    the WENO reconstructions (it can also generate C and OpenCL
    kernels).  The generated kernels
    can optionally perform individual steps of the reconstruction
    process (i.e., computing smoothness indicators, nonlinear
    weights, or the final reconstruction) or combinations thereof.
    This affords authors some flexibility in avoiding redundant
    computations or minimizing memory allocations and accesses.  Furthermore, the
    points within each cell at which the WENO reconstruction is
    performed are also arbitrary.  Negative weights are automatically split
    into positive and negative parts \cite{Shi2002}, allowing PyWENO to generate kernels for
    routines for arbitrary sets of points (such as arbitrary order
    Gauss-Legendre quadrature points).

    For PyClaw, odd-order WENO routines to approximate the
    solution at the left and right edges of each cell were generated
    from fifth to seventeenth order. All aspects of the WENO
    reconstruction are wrapped into standalone subroutines and no
    temporary work arrays are allocated.  Using these routines is
    trivially easy; for instance, to use the 9th-order WENO method
    instead of the classic algorithm in Listing \ref{hello-world} one needs only to replace line
    9 by the two lines
    
    \begin{lstlisting}[caption=Using SharpClaw (with PyWENO)]
    solver = pyclaw.SharpClawSolver1D()
    solver.weno_order = 9
    \end{lstlisting}

\section{Parallelization\label{sec:petclaw}}
Like many finite difference and finite volume codes, \clawpack and SharpClaw
implement boundary conditions through the use of {\em ghost cells}.  In this
approach, fictitious layers of cells are added around the edge of the problem
domain; the number of layers depends on the width of the stencil of the
numerical scheme.  At the beginning of each step, the ghost cell values are set
to satisfy the specified boundary conditions.
Then the numerical scheme is applied on all the interior (non-ghost) cells.
Many types of
boundary conditions can be handled well in this manner, including periodicity,
reflection, and outflow (non-reflecting).  Custom boundary conditions may also
be specified, for instance to model time-dependent inflow.

This approach is highly amenable to parallelization, since it is based on
the idea that information at the edge of a domain is filled in by a routine that is
independent of the rest of the numerical scheme. Therefore, the serial kernels
can be applied on each processor of a distributed
parallel machine as long as some routine first fills the ghost cells on the processor
either by appeal to boundary conditions or through communication with neighboring processors,
as appropriate. Only this ghost cell routine needs to know the global topology of
the problem; the serial kernels operate based entirely on local information.
This orthogonality allows independent development of serial numerical schemes
and parallel communication patterns, and is a key strategy in combining the
work of computational mathematicians and computer scientists.

The same global-local decomposition is employed in PETSc.
The PETSc library includes a DMDA object that implements parallelization through the use of ghost cells. The DMDA is
a highly scalable class for data layout across parallel, structured grids. All storage allocation and communication of
ghost values is handled by the DMDA, but storage is returned to the PyClaw program as numpy arrays so that no code
changes are necessary and the user experience is identical. In Figure \ref{fig:DAnumbering}, we show three different
representations of data over a simple $5 \times 6$ structured grid. The global ordering is used as input to PETSc linear
and nonlinear solvers, whereas the natural ordering is used for output since it is independent of the particular
parallel partition. Local ordering is used to extract data over a ``halo'' region, including ghost unknowns shared with
other processes.

This is, in fact, how PyClaw makes use of the DMDA structure. Local vectors are extracted with a given number of overlap
unknowns, and computations are performed using the same serial routines. These local vectors are then used to update a
global vector, and PETSc performs the appropriate accumulation for shared unknowns. This simple mechanism in PETSc for
integrating local and global data (which works also for unstructured grids) allows easy parallelization.
Thus PyClaw relies on \clawpack and \sharpclaw
to provide computational kernels for time-dependent nonlinear wave propagation
and on PETSc (through petsc4py) to manage distributed data arrays and the communication between them.
The data structures in PETSc and \clawpack/\sharpclaw are directly interfaced through the
Python package numpy \cite{numpy}. 

The parallel extension of PyClaw consists of only about 300 lines of Python code. 
Any PyClaw script can be run in parallel simply by replacing the statement \texttt{from clawpack import pyclaw} with 
\begin{lstlisting}[caption=Running in parallel]
from clawpack import petclaw as pyclaw
\end{lstlisting}
and invoking the Python script with \texttt{mpirun}.


\begin{figure}
  \centering
\setlength{\tabcolsep}{5pt}
\begin{tabular}{ccc}
\begin{tabular}{c}
\begin{tabular}{|ccc|cc|}
\hline
\multicolumn{3}{|c|}{Proc 2} & \multicolumn{2}{c|}{Proc 3} \\
\hline
25 & 26 & 27 & 28 & 29 \\
20 & 21 & 22 & 23 & 24 \\
15 & 16 & 17 & 18 & 19 \\
\hline
10 & 11 & 12 & 13 & 14 \\
 5 &  6 &  7 &  8 &  9 \\
 0 &  1 &  2 &  3 &  4 \\
\hline
\multicolumn{3}{|c|}{Proc 0} & \multicolumn{2}{c|}{Proc 1} \\
\hline
\end{tabular} \\
Natural numbering
\end{tabular}
&
\begin{tabular}{c}
\begin{tabular}{|ccc|cc|}
\hline
\multicolumn{3}{|c|}{Proc 2} & \multicolumn{2}{c|}{Proc 3} \\
\hline
21 & 22 & 23 & 28 & 29 \\
18 & 19 & 20 & 26 & 27 \\
15 & 16 & 17 & 24 & 25 \\
\hline
 6 &  7 &  8 & 13 & 14 \\
 3 &  4 &  5 & 11 & 12 \\
 0 &  1 &  2 &  9 & 10 \\
\hline
\multicolumn{3}{|c|}{Proc 0} & \multicolumn{2}{c|}{Proc 1} \\
\hline
\end{tabular}\\
Global numbering
\end{tabular}
&
\begin{tabular}{c}
\begin{tabular}{|ccc|cc|}
\hline
\multicolumn{3}{|c|}{Proc 2} & \multicolumn{2}{c|}{Proc 3} \\
\hline
 X &  X &  X &  X &  X \\
 X &  X &  X &  X &  X \\
\gray{12} & \gray{13} & \gray{14} & \gray{15} &  X \\
\hline
 8 &  9 & 10 & \gray{11} &  X \\
 4 &  5 &  6 &  \gray{7} &  X \\
 0 &  1 &  2 &  \gray{3} &  X \\
\hline
\multicolumn{3}{|c|}{Proc 0} & \multicolumn{2}{c|}{Proc 1} \\
\hline
\end{tabular} \\
Local numbering
\end{tabular}
\end{tabular}
  \caption{A simple $5 \times 6$ structured mesh is represented using a DMDA. The leftmost figure shows the natural
    numbering of unknowns, for instance used by most visualization packages, which is independent of the parallel
    layout. The middle figure shows the PETSc global numbering of unknowns, which is contiguous for each process. The
    rightmost figure shows the local numbering of unknowns for process 0, which includes \gray{ghost} unknowns shared
    with neighboring processes.}
  \label{fig:DAnumbering}
\end{figure}

The serial PyClaw routines handle discretization, Riemann solves, limiting and reconstruction, since they only depend on
local data. PETSc handles parallel layout and communication, but has no information about the
local computations. PETSc allows fine-grained control of the ghost value communication patterns so that parallel
performance can be tuned to different supercomputing architectures, but by default a user does not need to manage parallelism or see PETSc code. In fact, the PetClaw user is
shielded from PETSc in much the same way that a PETSc user is shielded from MPI. This separation can enable future
development. For instance, an unstructured mesh topology of hexahedral elements could be managed by PETSc, using a Riemann
solver which could accommodate deformed elements, without changes to PyClaw.

In addition to communication of ghost cell values, parallel hyperbolic solvers require
communication of the maximum wave speed occurring on each processor in order to
check whether a prescribed stability condition (generally phrased in terms of the Courant number)
has been satisfied and choose the size of the next time step appropriately.
This is also handled by a single PETSc call.

Although the basic \clawpack and SharpClaw algorithms are explicit and
require no algebraic solver, a powerful advantage gained by using PETSc for
parallelization is the possibility of employing PETSc's solvers for implicit
integration of hyperbolic problems that are stiff due to source terms or fast waves that
do not need to be accurately resolved.  This is the subject of ongoing work.

A particular challenge of using Python is that most parallel debuggers support
only C or Fortran, making parallel debugging of Python codes difficult \cite{Enkovaara2011}.
This is yet another motivation for using a tested parallel library like PETSc.

\section{Performance\label{sec:performance}}
    A few previous works have considered efficiency of scientific Python
    codes in serial as well as in parallel;
    see for instance \cite{Cai2005,Langtangen2008,Nilsen2010,Enkovaara2011}.
    Those studies consisted
    mainly of simple code snippets run in serial
    or on up to a few dozen processors until the recent work
    \cite{Enkovaara2011}, which includes scalability
    studies up to 16,384 cores.  In this section, we investigate
    the efficiency of a full object-oriented Python framework (PyClaw)
    compared with hand-coded Fortran (\clawpack).  We also consider the scaling of
    PetClaw on all 65,536 cores of the Shaheen supercomputer at KAUST.

    We consider only the second-order classic Clawpack algorithm here, as we are mainly
    interested in the effect of using a Python framework (in the serial
    case) and the cost of communication (in the parallel case). In terms
    of these factors, roughly similar results may be expected for the
    performance of the higher order algorithms, and preliminary tests
    (not described here) indicate good scaling of those also. 

    \subsection{Serial performance}
    For a detailed serial performance comparison of an explicit stencil-based
    PDE code in Python, see \cite{Langtangen2008}. In that work, vectorized numpy
    code was found to be fast enough for some operations, while wrapped Fortran
    loops performed identically to a pure Fortran code.  In contrast to the
    simple kernel code considered there, we present tests of a full object-oriented
    solver framework.  Our results thus extend those of \cite{Langtangen2008},
    providing an indication of the efficiency that can be expected for a
    sophisticated Python-based PDE solver framework.

    Table \ref{2DSerialClawpackPetclawComparison} shows an on-core serial
    comparison between
    the Fortran-only \clawpack code and the corresponding hybrid PyClaw
    implementation for two systems of equations on two different platforms.
    The hyperbolic systems considered are the 2D linear acoustics equation and the 2D
    shallow water (SW) equations \cite{levequefvmbook}. The acoustics test involves a very simple
    Riemann solver (amounting to a $3\times 3$ matrix-vector multiply)
    and is intended to provide an upper bound on the performance loss
    arising from the Python code overhead. The shallow water test involves a
    more typical, costly Riemann solver (specifically, a Roe solver with an
    entropy fix) and should be considered as more representative of realistic
    nonlinear application problems.
    Clawpack and PyClaw rely on similar Fortran kernels that differ only in the array
    layout.  Because most of the computational cost is in executing the
    low-level Fortran kernels, the difference in performance is relatively
    small -- though not negligible.  The results for the Shallow water equations are
    in rough agreement with the 10\% overhead reported in \cite{Enkovaara2011}.
    A 10-30\% increase in computational time (for realistic applications)
    seems well worth the advantages provided by the use of Python (in
    particular, easy parallelism).  The overhead is expected to be even smaller
    for more complex systems or in three dimensions.

    \begin{table}[h]
    \center
    \caption{Timing results (in seconds) for on-core serial experiments solving
    acoustics and shallow water problems implemented in both \clawpack and
     PyClaw on Intel Xeon and the IBM BlueGene/P PowerPC 450 processors}
    \label{2DSerialClawpackPetclawComparison}
    \begin{tabular}{|c|c|c|c|c|}
    \hline
    \textbf{Application} & \textbf{Processor}  & \textbf{\clawpack} & \textbf{PyClaw} & \textbf{Ratio} \\ \hline
    \multirow{2}{*}{Acoustics}     & Intel Xeon  & 28s  & 41s  & 1.5 \\
                                   & PowerPC 450 & 192s & 316s & 1.6 \\ \hline
    \multirow{2}{*}{Shallow Water} & Intel Xeon  & 79s  & 99s  & 1.3 \\
                                   & PowerPC 450 & 714s & 800s & 1.1 \\ \hline
    \end{tabular}
    \end{table}

    \subsection{Parallel performance}
    We now investigate the parallel performance of PetClaw on the Shaheen supercomputer at KAUST,
    an IBM BlueGene/P system consisting of 65,536 cores.  When characterizing the performance of scientific codes on
    supercomputers, a commonly used characterization is that of \textit{weak scalability}, which is assessed by studying
    how the  run time of the code is affected when the resolution of the simulation and the number of processors is
    increased commensurately to maintain a fixed amount of work per processor.  The \texttt{parallel efficiency} is given by
    dividing the run time of the single processor job by the run time of the parallel job.  

\begin{figure} \centering
        \includegraphics{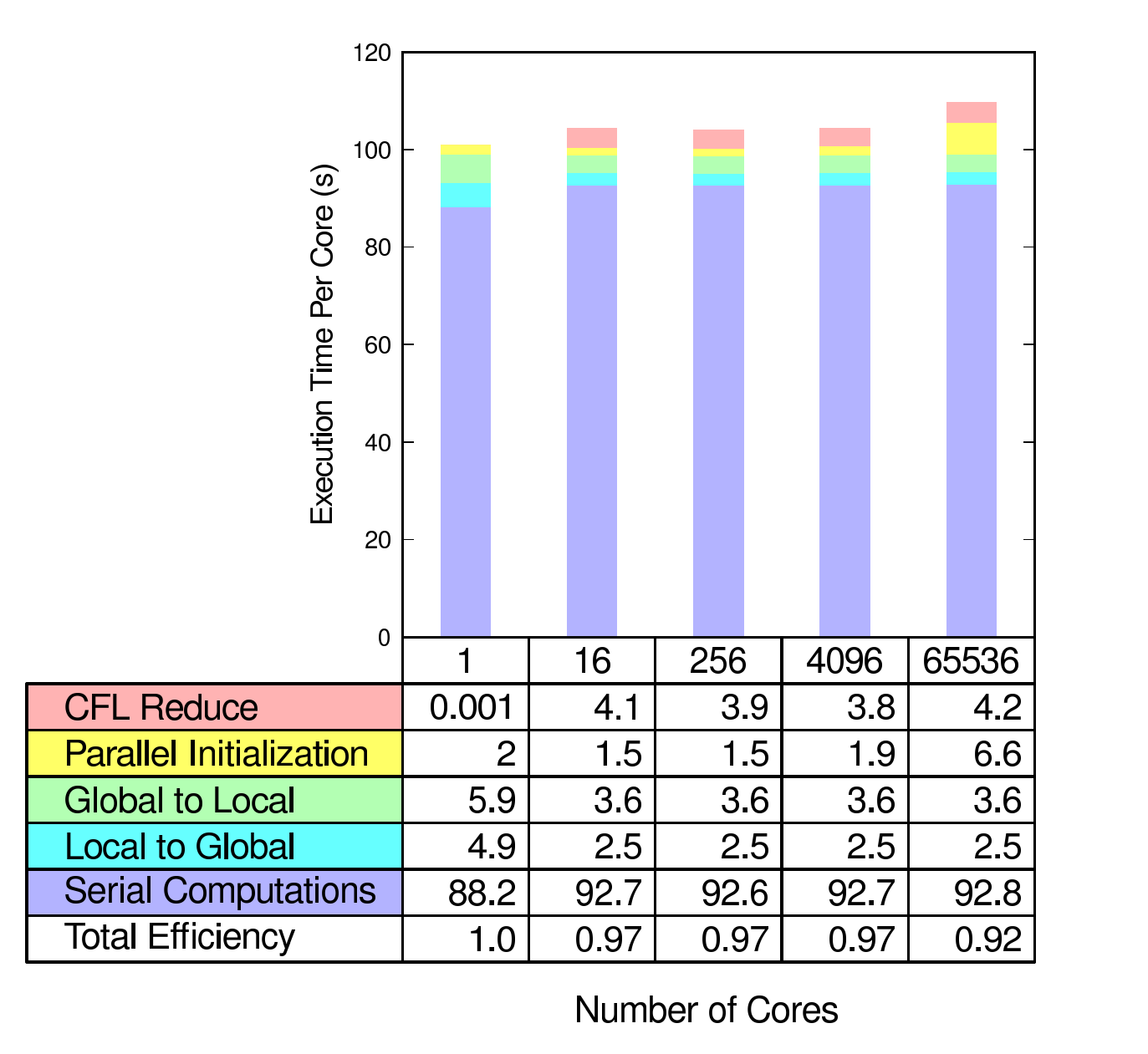}
        \caption{Weak scaling performance profile of the shock bubble problem with 160,000 grid cells per core\label{fig:par_prof_euler}}
    \end{figure}

    The problem
    used for the comparisons is of a compressible, inviscid flow that consists of a shock
    impacting a low-density bubble (examined in detail in Section \ref{sec:shockbubble}).
    We investigate weak scaling by running the problem for a fixed number of time steps and with a fixed
    number of grid cells ($400\times 400=160,000$) per core, while increasing the number of
    cores from one up to the whole machine. Figure \ref{fig:par_prof_euler}
    shows the results, with parallel efficiency provided in the last row.
    It is important to note that the time required to load the necessary
    Python packages and shared objects, which occurs only once at the
    beginning of a simulation (or series of batch simulations) has been
    excluded from the results presented here.  This load time is
    discussed in the next section.

    Observe that in all parallel
    runs, more than 90\% of the time is spent in the computational kernels.
    The parallel operations scale extremely well: the the CFL condition-related reduction
    takes essentially the same amount of time for all runs from 16 processors
    up, as does the communication of ghost cell values in localToGlobal.
    Together these parallel operations consume about 6\% of the total
    run time.   Parallel initialization consists of PETSc parallel object
    construction, including memory allocation and MPI communicator
    initialization. Note that the parallel initialization, while significant in
    these artificial test runs, will not contribute significantly to the cost of
    real simulations because it is a one-time cost.

    \subsection{Dynamic Loading}

    As alluded to already, the work of loading Python libraries and
    objects dynamically at run-time does not currently scale well on the Shaheen
    system.  Large-scale supercomputers such as Shaheen rely on parallel file
    systems that are designed to support large distributed loads, with each
    process independently accessing data.  Dynamic loading does not follow this
    pattern because every process is attempting to access the same data
    simultaneously.  This issue was partially addressed in
    \cite{Enkovaara2011}, but an implementation capable of
    supporting dynamic library loading is still lacking.

    The dynamic loading time for the PetClaw runs in Section \ref{sec:apps}
    is less than 5\% of the total simulation time, and this will generally
    be the case for 2D wave propagation problems because the CFL condition
    means that large simulations of
    hyperbolic problems necessarily require long run times in order for waves
    to propagate across the full domain.


\section{Software Development Practices\label{sec:software}}
    The majority of software development practices utilized in PyClaw are
    inherited from the open source software community.  The community's
    atmosphere of free and open sharing complements the tradition of scientific inquiry.  In
    fact, growing movements within the scientific community seek to embrace \textit{scientific
      reproducibility} for software tools used in conducting mathematical and scientific research~\cite{fomel2009,Knepley2010}.

    In addition to making our results reproducible, we also intend that our software be
    useful as a reference for understanding numerical methods involved in solving hyperbolic PDEs and as a platform for extending and applying these techniques.  As such, we also seek to provide a
    friendly and inviting context for scientists working in this cross-disciplinary environment to conduct their
    research. 

  \subsection{Application Tests}
  The goal of reproducibility in research is to improve not only confidence in results, but their
  extensibility as well.  A series of regression tests have been devised for every major application
  where PyClaw has been used.
  The script and parameters for generating the test are stored in the repository (typically in
  a single file), along with verified output from a successful run.  Where possible, the output is
  verified against a different solver implementation or by analysis.  These ``application tests'' produce the
  same output regardless of the choice of solver type, kernel implementation, or computer architecture.  The python-nose
  (\url{http://code.google.com/p/python-nose/}) unit testing framework simplifies development and
  selective execution of the tests.  Programmatic test code generation is used to exercise the full range of solver and
  kernel options for each test.  Scientific results are archived as application tests within the
  unit testing framework, ensuring that our published results are reproducible in current and future versions of the
  PyClaw solver.

  In our experience, the application tests are the single greatest factor in facilitating
  adoption, refactoring, and extension of the code.  New users are confident that they have a working installation (or
  can tell us what doesn't work on their architectures) and are capable of reproducing our published scientific results.
  Developers refactoring the code for improved organization, performance, or readability can rely on the application
  tests for regression testing, to ensure that their changes have not incidentally  broken anything.  Perhaps most
  importantly, new solver methods and their implementations can be verified against known solutions with the application
  tests.  This facilitates and encourages the development of new ideas within the PyClaw framework.

  \subsection{Hosted Distributed Version Control}
  Our use of git (\url{http://git-scm.com/}), a modern, distributed version control system, provides many benefits. Development need not be
  synchronized through a master server, which makes it easier to incorporate subprojects from developers loosely
  attached to the core team. Management of releases and bugfix updates has been greatly simplified. However, perhaps the
  greatest beneficiary is the user. Users do not have to wait for PyClaw releases in order to retrieve bugfixes for
  particular machines or improvements which are under active development, they need only update to a given
  changeset. Moreover, a user can easily switch between an approved release and experimental code for comparison with a
  single version control command. This allows the user a much finer-grained manipulation of versioning than was
  previously possible.

  There are many excellent open source distributed version control hosting sites, including Bitbucket
  (\url{http://www.bitbucket.org}) and GitHub
  (\url{http://www.github.org}), which provide a range a services to both developers and
  community users. PyClaw leverages the services provided at GitHub, which includes
  wiki webpages for user
  communication, as well as a bug reporting and tracking infrastructure integrated with the hosted version control repository.
  We have separately engaged the use of Google Groups to provide a mailing list for the
  PyClaw user and developer community.

 \subsection{Documentation}
  PyClaw is provided with a range of documentation suitable for the variety of users interacting with the software.
  While this paper provides a high-level overview of the capabilities of the code and its application, it is our
  experience from using other projects that the best software documentation includes a separate tutorial and user's
  guide with a class reference section.  The tutorial and user's guide are maintained in the ReStructured Text format, from which they
  can be translated into HTML, PDF, and several other output formats using for instance Sphinx (\url{http://sphinx.pocoo.org/}).  The PyClaw code itself is documented inline using Python's docstring
  conventions, allowing us to automatically generate class and module reference sections for our documentation.

\section{Applications\label{sec:apps}}
    The numerical algorithms made accessible in PyClaw, empowered by parallelization,
    are capable of modeling challenging wave propagation phenomena.  In this section,
    we provide three example applications.  They are not intended to break new scientific
    ground, but rather to demonstrate the versatility
    of the algorithms accessible through PyClaw, and (in the third application) the
    power of PetClaw as a scalable parallel solver.

    \subsection{Shallow Water Flow on the Sphere}\label{sub:shallow_sphere}
     Classical shallow water equations on a sphere are an approximation of the
     flow on the earth's surface. They are of interest because they capture most
     of the flow's features of a thin layer of fluid tangent to the surface of
     the sphere. Therefore, they are often used by meteorologists,
     climatologists and geophysicists to model both atmosphere and oceans.

     In three-dimensional Cartesian coordinates, using $h$ and
     $\mathbf{u} = \left(u,v,w\right)^T$ to define the height and the fluid velocity vector,
     the shallow water equations are of the form \eqref{eq:conslaw}, with
     $\mathbf{q} = \left(h, hu, hv, hw\right)^T$
     and
     \begin{equation}
        \mathbf{f}\left(\mathbf{q}\right) =
            \left(\begin{array}{ccc}
                    hu & hv & hw \\
                    hu^2 + \frac{1}{2}gh & huv & huw \\
                    huv  & hv^2 + \frac{1}{2}gh & hvw \\
                    huw  & hvw  & hw^2 + \frac{1}{2}gh
                   \end{array}\right),
     \end{equation}
     where $g$ is the gravitational acceleration. The source term
     $\mathbf{s}\left(\mathbf{q},\mathbf{x}\right)$ includes the
     Coriolis force and an additional term
     that ensures that the velocity is tangent to the sphere:
     \begin{equation}\label{eqn:sourcesweqns}
     \mathbf{s}\left(\mathbf{q},\mathbf{x}\right) =
     -\frac{2 \Omega}{a} z \left(\mathbf{x} \times h\mathbf{u}\right) +
     \left(\mathbf{x} \cdot \left(\mathbf{\nabla} \cdot \tilde{\mathbf{f}}
     \right)\right) \mathbf{x}.
     \end{equation}
     Here $\Omega$ and $a$ are the rotation rate and the radius of the earth,
     respectively. In \eqref{eqn:sourcesweqns} $\tilde{\mathbf{f}}$ is the part
     of the flux matrix associated with the momentum equations
     \cite{D.A.Calhoun2008}.

     In this framework, we consider the numerical solution of the zonal wave
     number 4 Rossby-Haurwitz problem \cite{Williamson1992}. 
     Rossby-Haurwitz waves are steadily propagating initial conditions of the 
     nonlinear non-divergent barotropic vorticity equation on a rotating sphere
     \cite{Haurwitz1940}.   Although they do not represent exact steady solutions
     of the shallow water equations, they are expected to evolve in a similar way
    \cite{Thuburn2000}. 
     For this reason Rossby-Haurwitz waves have also been used to test shallow 
     water numerical models and are among the standard shallow water model test 
     cases proposed by Williamson et al. \cite{Williamson1992}.
       
     The problem consists of
     a divergence-free initial velocity field that rotates about the
     $z-$axis without significantly changing form on short time scales (dynamical 
     weak instabilities). On 
     longer time scales (50-100 days) the instabilities effects lead to the 
     breakdown of the wave structure. Previous studies have shown that the time 
     at which the solution symmetry breaks depends strongly on the numerical
     truncation error of the scheme \cite{Thuburn2000}. Because of this
     the Rossby-Haurwitz problem is frequently used to assess
     the accuracy of a numerical algorithm for the solution of the shallow water
     equations on a rotating sphere.

     The three-dimensional shallow water equations are solved on the logically rectangular
     grid introduced in \cite{D.A.Calhoun2008}, using the same approach employed there,
     namely the classic Clawpack method with Strang splitting for the source term \cite{G.1968}.
     Simulations have been performed on four grids with $100 \times 50$,
     $200 \times 100$, $400 \times 200$ and $800 \times 400$ cells. Table
     \ref{tbl:rossby-haurwitz} lists the breakdown time. With the coarsest mesh,
     the diffusive part of the numerical error suppresses the instability completely
     (cfr. \cite{Thuburn2000}).
     The finer grid results confirm that the time at which the initial velocity loses its
     symmetry is sensitive to the numerical error.
     \begin{table}[h]
     \center
     \caption{Time at which the symmetry of the
     Rossby-Haurwitz wave breaks down. \label{tbl:rossby-haurwitz}}
     \begin{tabular}{c|c}
     \hline
     Grid & Breakdown \\
     \hline
     $100 \times 50$   &            - \\
     $200 \times 100$  & $\approx$ 34 d \\
     $400 \times 200$  & $\approx$ 45 d \\
     $800 \times 400$  & $\approx$ 46 d \\
     \hline
     \end{tabular}
     \end{table}

\begin{figure}
    \centering
    \subfigure[0 days.]{\includegraphics[scale=0.12]{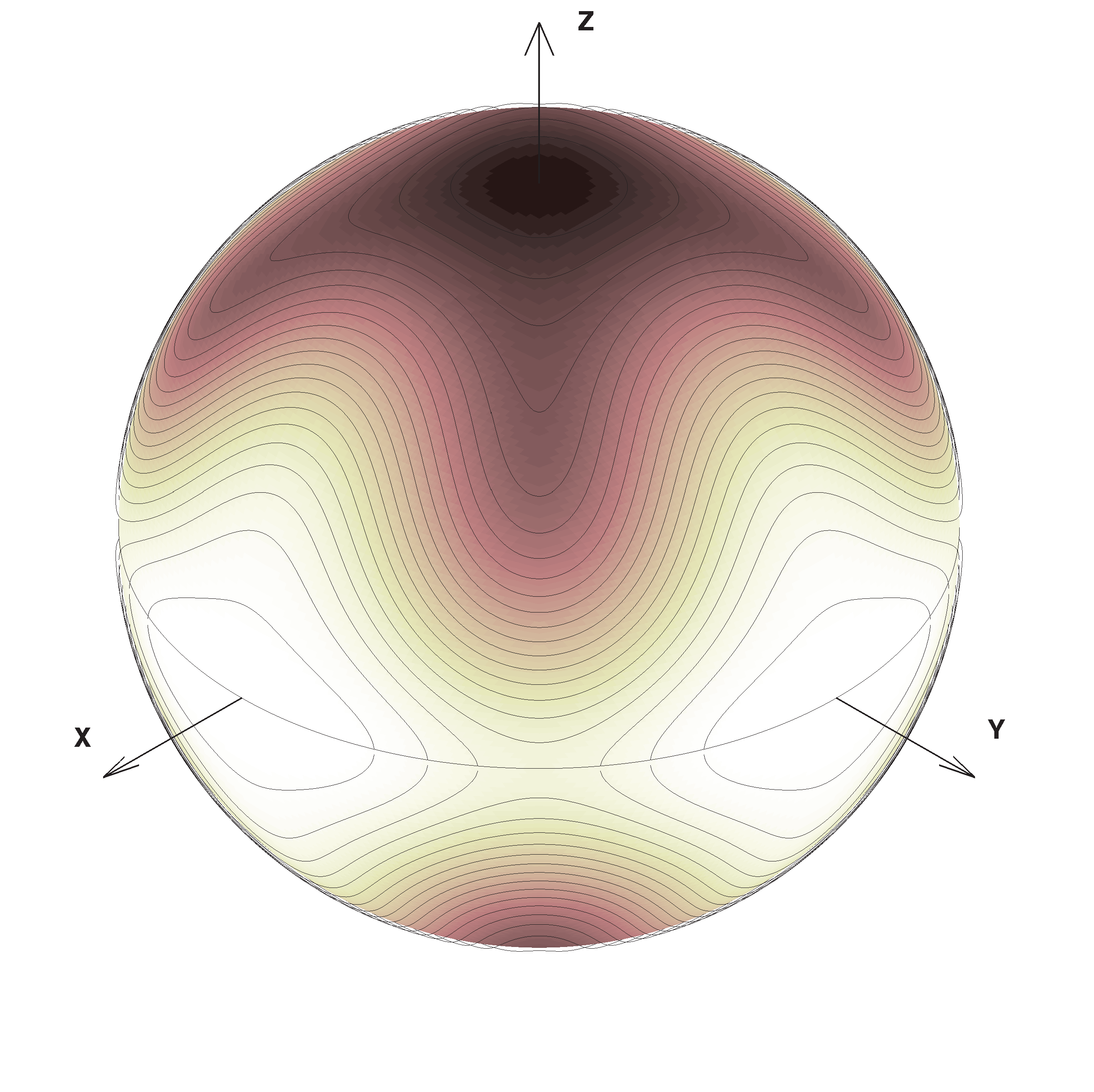}} \quad
    \subfigure[38 days.]{\includegraphics[scale=0.12]{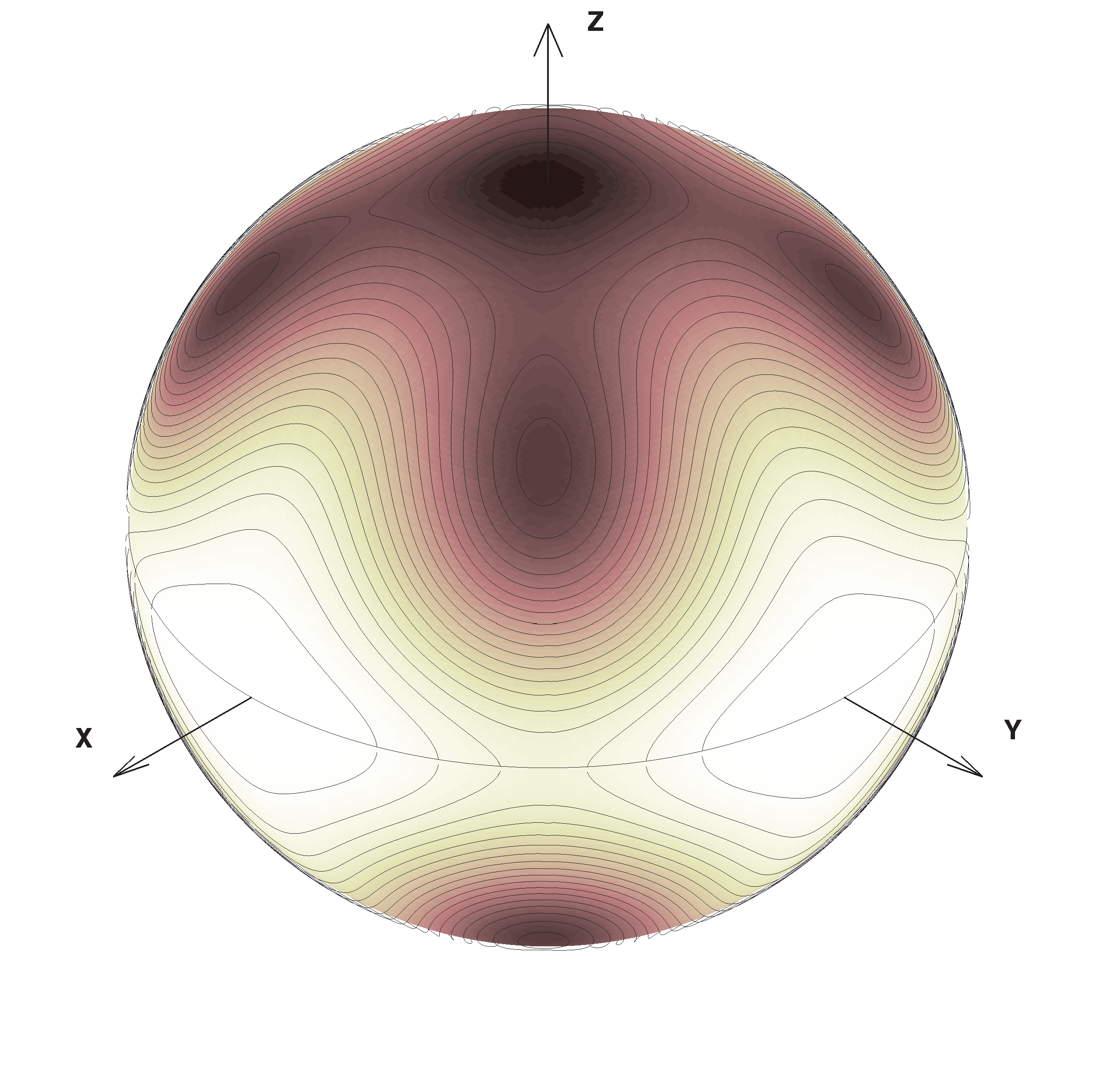}}
    \subfigure[45 days.]{\includegraphics[scale=0.12]{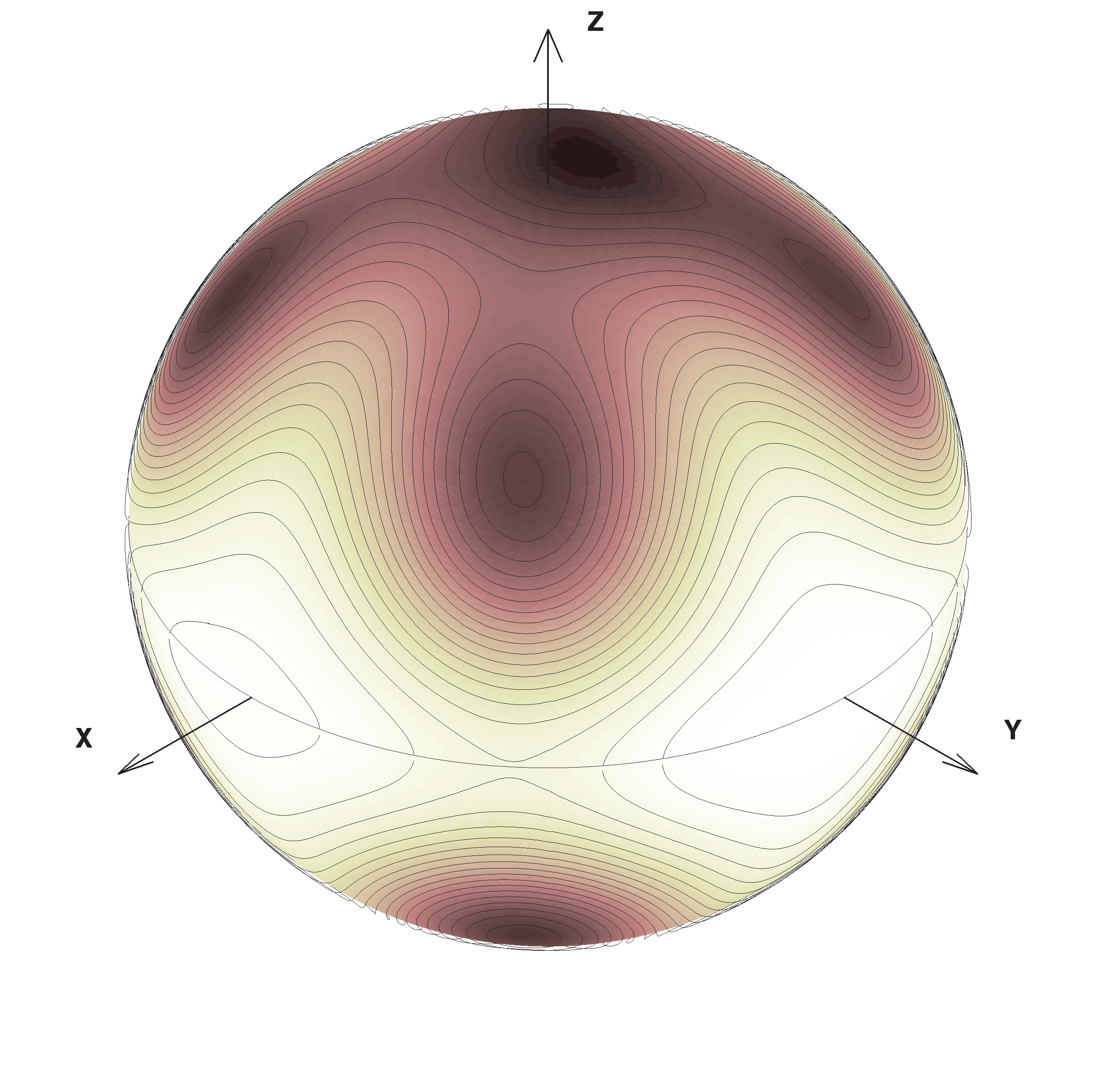}} \quad
    \subfigure[48 days.]{\includegraphics[scale=0.12]{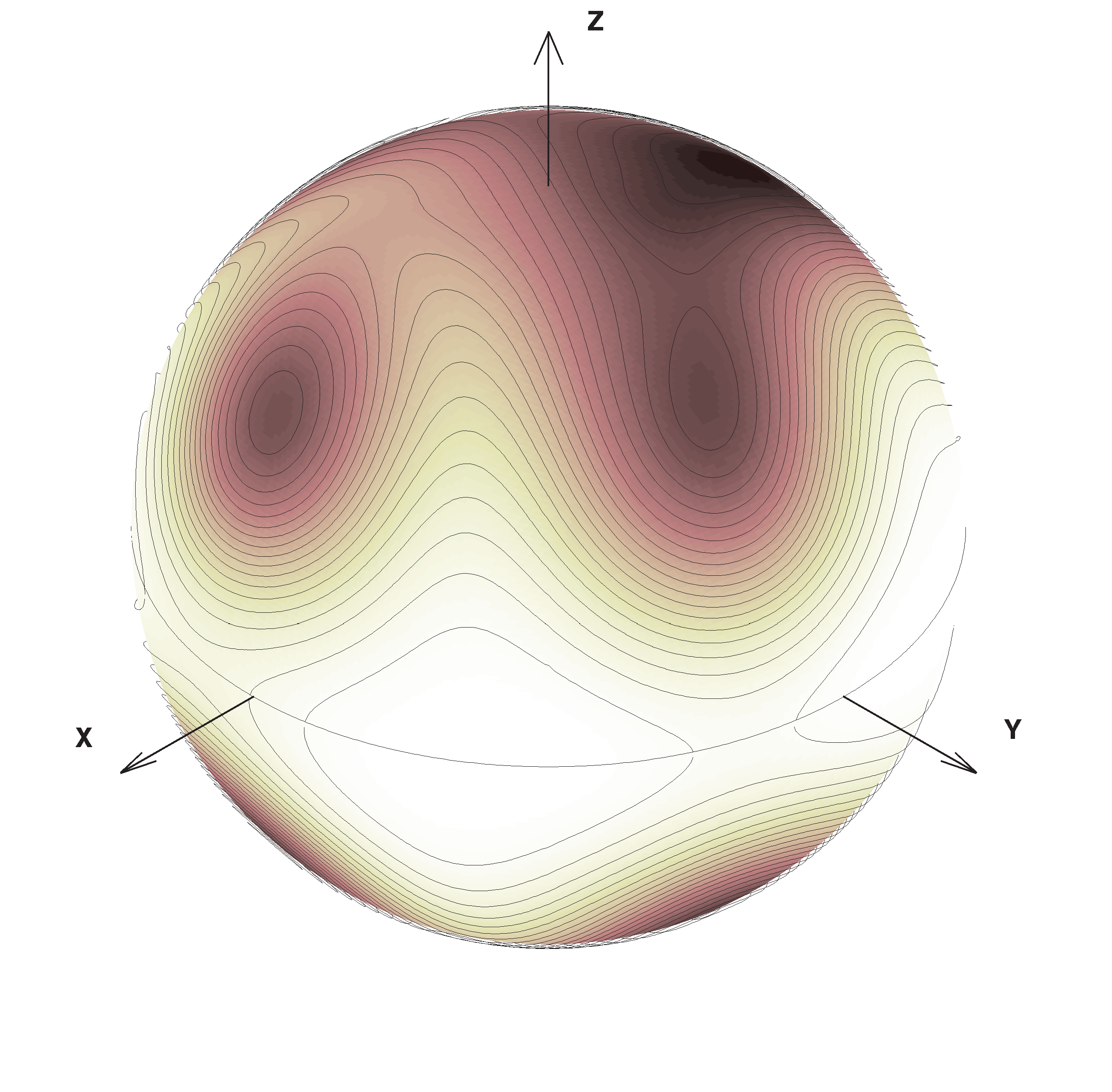}}
    \caption{Water height of the zonal wave number 4 Rossby-Haurwitz problem at
    different times; grid with $400 \times 200$ cells. Contour interval 120 m.
    \label{fig:SWSphereRossby-Haurwitz}}
    \end{figure}

     Figure \ref{fig:SWSphereRossby-Haurwitz} shows the contour line of the
     water height at day 0 (initial solution), 38, 45 and 48, for the grid with $400 \times 200$
     cells. These plots show qualitatively the evolution of the instability.

    \subsection{Shock-Bubble Interaction\label{sec:shockbubble}}
    The Euler equations for a compressible, inviscid fluid with cylindrical symmetry
    can be written as
    \begin{equation} \label{euler-rz}
    \begin{aligned}
        \rho_t + (\rho u)_z + (\rho v)_r & = -\frac{\rho v}{r}, \\
        (\rho u)_t + (\rho u^2 + p)_z + (\rho u v)_r & = -\frac{\rho u v}{r}, \\
        (\rho v)_t + (\rho u v)_z + (\rho v^2 + p)_r & = -\frac{\rho v^2}{r}, \\
        (\rho E)_t + \left((\rho E + p) u\right)_z + \left((\rho E + p) v\right)_r & =
                -\frac{(\rho E + p)v}{r}.
    \end{aligned}
    \end{equation}
    Here the $z$-coordinate represents distance along the axis of symmetry while
    the $r$-coordinate measures distance away from the axis of symmetry.
    The quantities $\rho, E, p$ represent density, total energy, and pressure, respectively,
    while $u$ and $v$ are the $z$- and $r$-components of velocity.

    We consider an ideal gas with $\gamma=1.4$ in the cylindrical domain
    $[0,2]\times[0,0.5]$.  The problem consists of a planar shock traveling in
    the $z$-direction that impacts a spherical bubble of lower-density fluid.
    In front of the shock $u=v=0$ and $\rho=p=1$ except inside the bubble, where $p=1,
    \rho=0.1$.  Behind the shock, $p=5,\rho\approx 2.82, v\approx 1.61$, and these
    conditions are also imposed at the left edge of the domain.
    In addition to \eqref{euler-rz}, we solve a
    simple advection equation for a tracer that is initially set to unity in
    the bubble and zero elsewhere in order to visualize where the fluid interior to the bubble
    is transported.  Reflecting boundary conditions are imposed at the bottom of the
    domain while outflow conditions are imposed at the top and right boundaries.

    Since the edge of the bubble is curved, it does not align with the Cartesian grid.
    Thus, in the cells that are partly inside and partly outside the bubble, the initial
    condition used is a weighted average of the different solution values, based on the
    fraction of the cell that is inside.  This fraction is computed by adaptive quadrature
    using the {\tt scipy.integrate.quad} package.

\begin{figure}
    \centering
    \subfigure[Initial tracer showing location of low-density bubble.]{\includegraphics[width=5in]{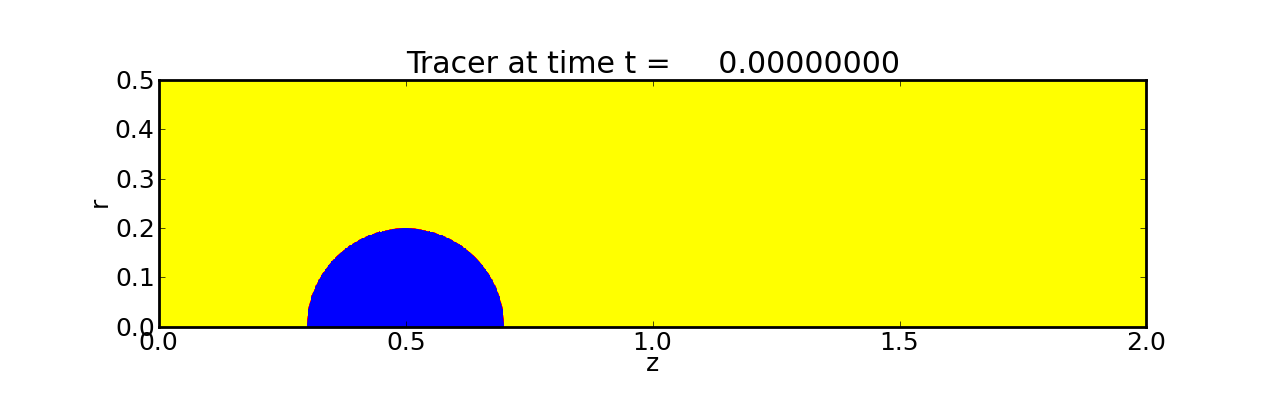}} \\
    \subfigure[Tracer showing location of of bubble material after shock impact.]{\includegraphics[width=5in]{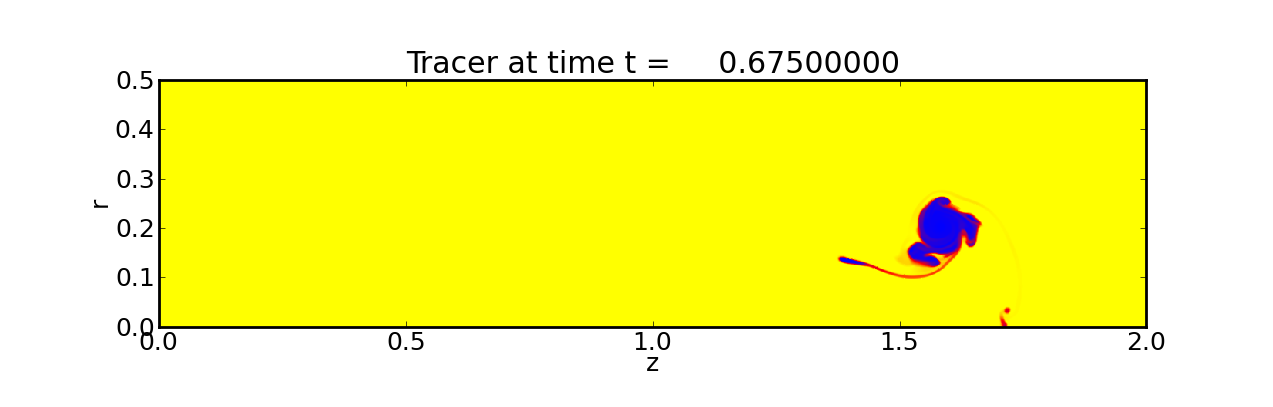}} \\
    \subfigure[Schlieren plot of density.]{\includegraphics[width=5in]{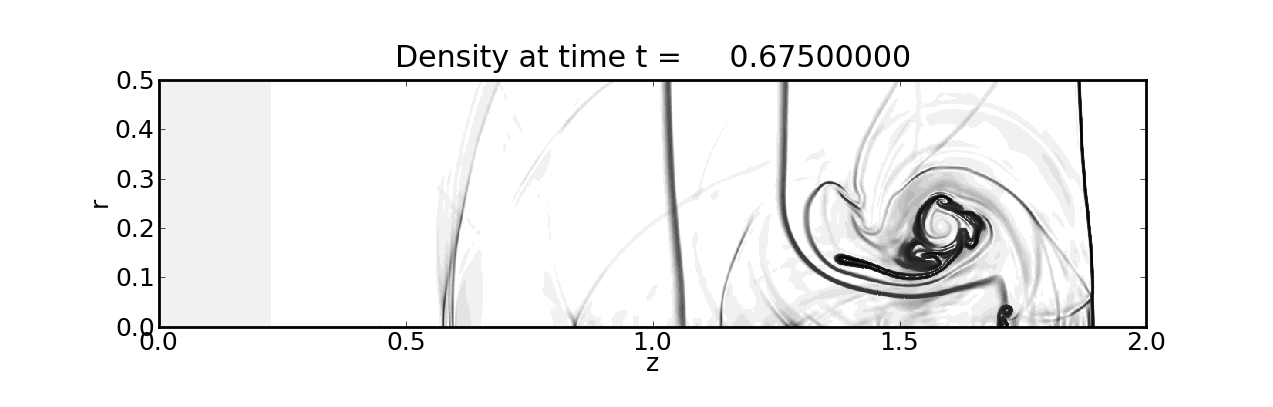}}
    \caption{Results of shock-bubble interaction computation, showing the transformation of
             the initially spherical bubble into a ``smoke ring'' after it is impacted by a shock
             wave.  \label{fig:shockbubble}}
    \end{figure}

\begin{figure}
    \centering
    \subfigure[Classic \clawpack solution.]{\includegraphics[width=2.5in]{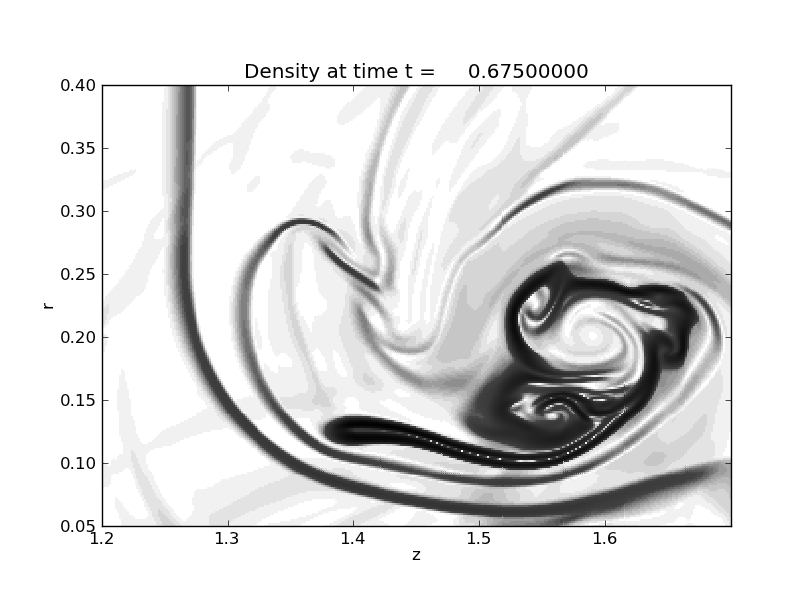}} \\
    \subfigure[\sharpclaw solution using WENO5 reconstruction.]{\includegraphics[width=2.5in]{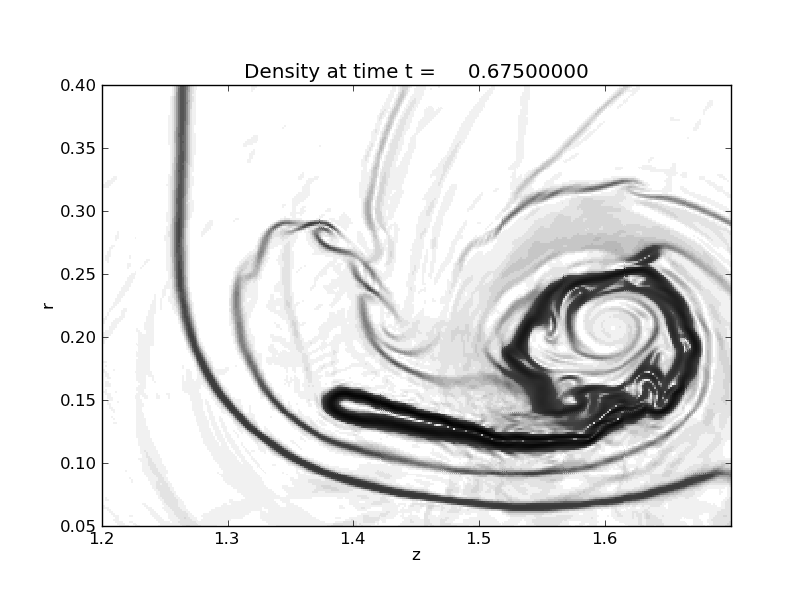}}
    \subfigure[\sharpclaw solution using WENO7 reconstruction.]{\includegraphics[width=2.5in]{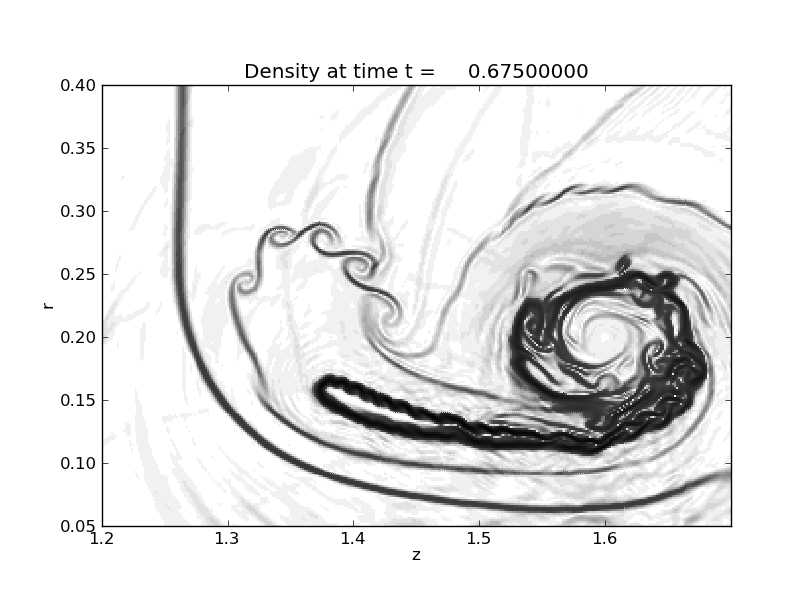}}

    \caption{Schlieren plots of density, zoomed in on the smoke ring.  All solutions are computed
            on a $1280\times 320$ grid. \label{fig:wenoshockbubble}}
    \end{figure}

    Figure \ref{fig:shockbubble} shows the initial condition and results of this problem,
    using a $1280\times 320$ grid and the unsplit classic Clawpack algorithm with full transverse corrections.
    The bubble is observed to transform into a ``smoke ring''.
    Considerably more detailed structure is evident in this simulation compared to
    lower-resolution adaptively refined results from AMRClaw that are published at
    \url{http://depts.washington.edu/clawpack/clawpack-4.3/applications/euler/2d/shockbubble/amr/www/index.html}.
    Figure~\ref{fig:wenoshockbubble} shows a close-up of the smoke ring solution obtained
    with the classic \clawpack algorithm, as well as solutions obtained
    using the \sharpclaw algorithm with fifth-
    (WENO5) and seventh-order (WENO7) reconstructions.
    All runs were performed on a $1280\times 320$
    grid with a maximum CFL number of 0.8.  Although the overall
    features of the solutions are similar, more fine structure is apparent
    in the \sharpclaw results.  For instance, several vortices can be seen
    to the left of the smoke ring in the WENO7 run that are not
    resolved in the classic run.

    \subsection{Cylindrical Solitary Waves in a Periodic Medium}\label{sec:periodic}
        The problem considered in this section is taken from \cite{QuezadadeLuna2011}.
        It involves the propagation of nonlinear waves in a two-dimensional crystal,
        leading to the formation of solitary radial waves or ``rings''.
        For this purpose, we consider the 2D $p$-system with spatially-varying
        coefficients as a model for elastic waves:
        \begin{equation} \label{p-system}
          \begin{aligned}
            \epsilon_{t}-u_{x}-v_{y} & = 0,\\
            \rho(x,y) u_{t}-\sigma(\epsilon,x,y)_{x} & = 0, \\
            \rho(x,y) v_{t}-\sigma(\epsilon,x,y)_{y} & = 0.
          \end{aligned}
        \end{equation}
        Here $\epsilon$ represents the strain, $u$ and $v$ are the velocities in $x$ and $y$
        respectively, $\rho(x,y)$ is the spatially-varying material density and
        $\sigma(\epsilon, x, y)$ is the stress. The system \eqref{p-system} is closed by
        introducing a nonlinear constitutive relation $\sigma(\epsilon, x, y)$. Similar to
        \cite{leveque2003}, we take
        \begin{equation}\label{nonlinear const relation}
          \sigma(\epsilon, x, y) = \exp(K(x, y)\epsilon)+1,
        \end{equation}
        where $K(x, y)$ is the spatially-varying bulk modulus.

        The medium, a small part of which is shown in Figure \ref{fig: checker
        medium}, is in a checkerboard pattern with alternating squares of two different materials:
        \begin{equation} \label{checker equation}
            (K(x,y),\rho(x,y)) = \left \{ \begin{aligned}
                    (1,1) & &\mbox{ if } \left(x-\lfloor x\rfloor-\frac{1}{2}\right)
                    \left(y-\lfloor y\rfloor-\frac{1}{2}\right)<0& \\
                    (5,5) & &\mbox{ if } \left(x-\lfloor x\rfloor-\frac{1}{2}\right)
                    \left(y-\lfloor y\rfloor-\frac{1}{2}\right)>0&.
            \end{aligned} \right .
        \end{equation}
         The problem is quite challenging, for multiple reasons.  First,
         the flux is spatially varying and even discontinuous -- meaning that
         the solution variables (strain and momentum) are also discontinuous.
         Furthermore, these discontinuities include corner singularities.
         Finally, and owing in part to these discontinuities, it is necessary to use
         a large number of cells per material period ($\gtrsim 100$)
         in order to get even qualitatively
         accurate results.  As we are interested in a phenomenon that arises
         only after waves pass through many ($>100$) material periods, this leads to
         a very high computational cost.

        The problem is solved using the \sharpclaw algorithm with fifth-order WENO
        reconstruction. As explained in Section \ref{sub:clawpack},
        the implementation in \sharpclaw is based on
        solving normal Riemann problems at the grid interfaces; see \cite{quezada2011}
        for a detailed explanation of the approximate Riemann solvers employed
        and for a much more detailed study of this problem.

\begin{figure}
          \begin{centering}
            \subfigure[Detail of checkerboard medium.]
              {\includegraphics[scale=0.3]{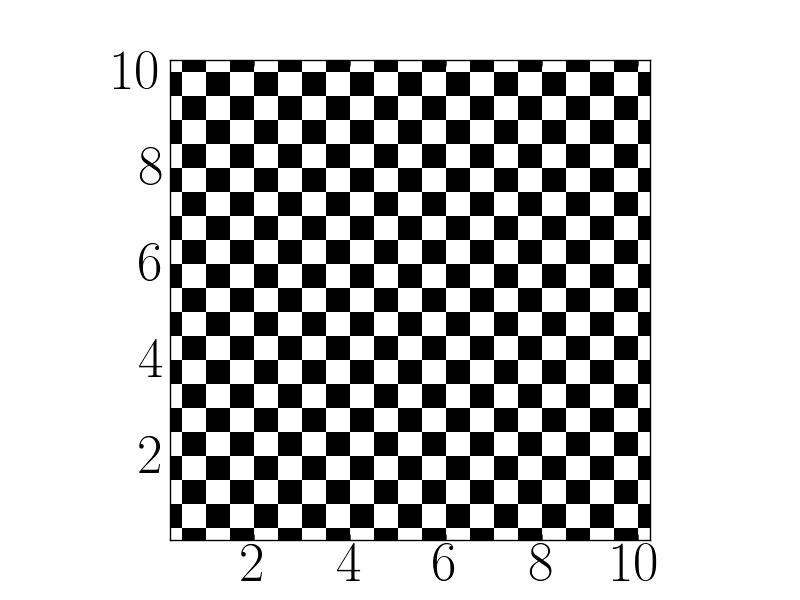}
                \label{fig: checker medium}}
            \subfigure[Detail of initial stress.]
              {\includegraphics[scale=0.3]{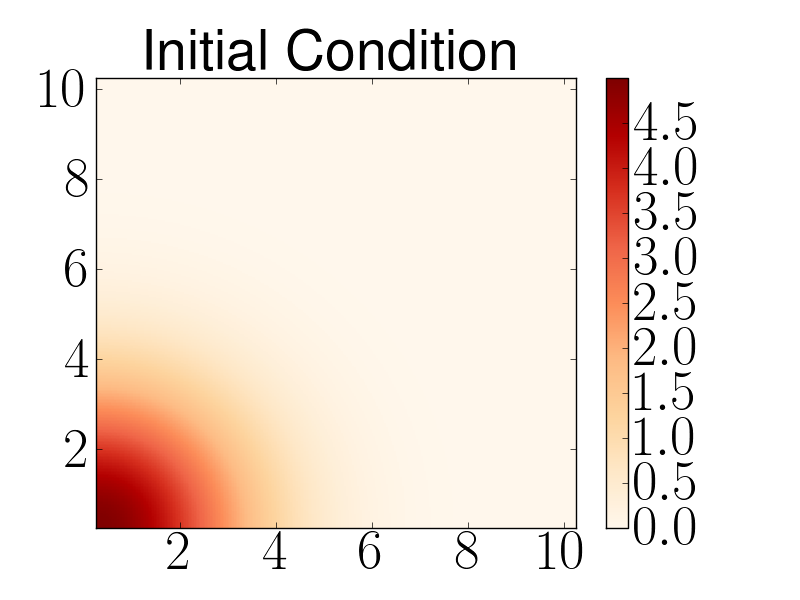}
                \label{fig: initial condition}}
          \par
          \end{centering}
          \caption{The medium and initial condition (both shown zoomed-in) for the cylindrical solitary wave problem.}
        \end{figure}

        A close-up view of the initial stress is shown in Figure~\ref{fig: initial condition}.
        The stress is a
        Gaussian pulse with an amplitude of $5$ and a variance in $x$ and $y$
        of $5$, centered at the origin.  The velocity is initially zero.
        The problem is symmetric with respect to reflection about the $x$- and $y$- axes,
        so the computational domain is restricted to the positive quadrant and reflecting boundary
        conditions are imposed at the left and bottom boundaries. Outflow (zero-order extrapolation)
        conditions are imposed at the top and right boundaries.  In units of the medium
        periodi, the domain
        considered is $200\times 200$, and the grid spacing is
        $\Dx=\Dy=1/240$.  Hence the full simulation involves $6.8\times 10^9$ unknowns.
        It was run on 16,384 cores of the Shaheen supercomputer at KAUST, over a period
        of about 3.2 days.

\begin{figure}
          \begin{centering}
            \subfigure[Solitary wave train (stress shown).]
              {\includegraphics[scale=0.3]{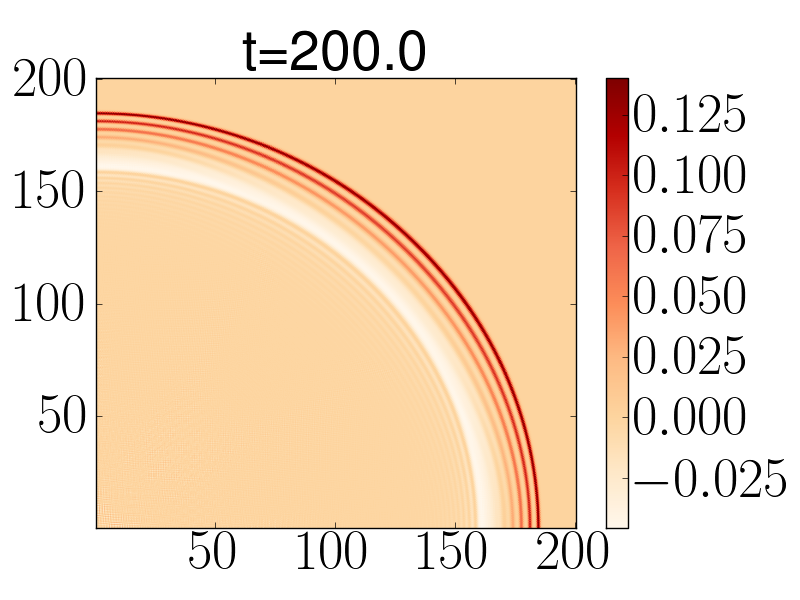}
                \label{fig: steg}}
            \subfigure[One-dimensional slices of the stress.  Solid black line: $45^\circ$;
                        dashed red line: $0^\circ$.]
              {\includegraphics[scale=0.3]{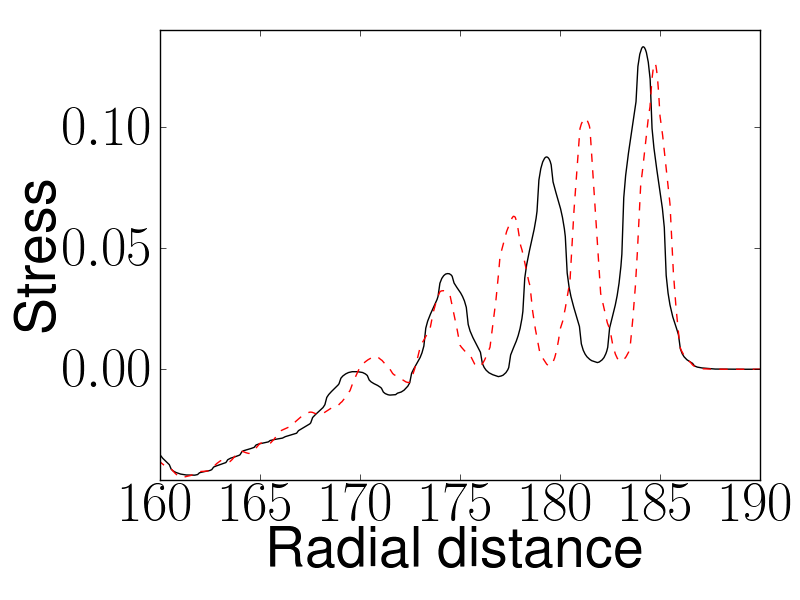}
                \label{fig: steg_slices}}
          \par
          \end{centering}
          \caption{Solution of the $p$-system \eqref{p-system} in a checkerboard medium,
          showing 2D cylindrical solitary wave formation.}
        \end{figure}

        The formation of solitary wave rings is seen clearly in Figure
        \ref{fig: steg}, which depicts the stress at $t=200$. The structure of
        these waves is shown in Figure \ref{fig: steg_slices}, which
        displays slices of the stress at $45^\circ$ (solid line) and $0^\circ$
        (dashed line) with respect to the $x$-axis.

\section{Discussion and future plans\label{sec:discuss}}
This work demonstrates that the use of Python
in combination with existing Fortran and C codes allows the
production of scientific software that is accessible, extensible, and efficient.
The serial performance loss is
relatively small, and is more than compensated for even on a typical laptop
by the ability to run in parallel without any additional effort.
Combining scalable parallelism with the algorithms of \clawpack and
SharpClaw yields a potent tool for exploring novel wave phenomena.

We are in the process of extending PyClaw to include
three-dimensional wave propagation and implicit time stepping.  Both
of these are straightforward steps, since the PyClaw framework
was written with such developments in mind and the related software
packages (\clawpack and PETSc) already support these features.
Preliminary implementations are under testing.

The use of Python in scientific computing has many potential
benefits \cite{perez2011}.  The Python packages numpy, scipy, and matplotlib offer essential
numerical tools with interfaces familiar to users of MATLAB (the {\em lingua
franca} of numerical methods) in a general-purpose programming language.
An increasing number of important libraries (like PETSc and Trilinos) now
have Python bindings, making it relatively easy to add powerful capabilities
like massively scalable parallelism to Python codes.  As discussed in Section
\ref{sec:software}, the
Python community promotes a range of positive code development practices
that are not common in scientific teams but are often picked up by those who
begin to work in Python \cite{FiPy:2009}.

While the existence of multiple scientific codes for solving the same problems
is healthy, it is widely recognized that increased code
sharing and reuse would benefit the numerical analysis and scientific
computing communities.  Closer integration of code developed by different
groups would not only allow researchers to be more productive (by reducing
duplication of effort), but would also allow useful algorithmic improvements to be
more rapidly distinguished from insignificant ones by simplifying the
task of comparing them.  In our experience, the adoption
of Python as a high-level scientific coding language dramatically increases
opportunities for code-sharing and reuse.  Indeed, the results described in
this paper consist largely of combining a few powerful existing pieces of
scientific software.

{\bf Acknowledgments.}  We are very grateful to Randall LeVeque,
and all the Clawpack developers, without whom PyClaw would not exist.
The authors also thank the following people:
\begin{itemize}
    \item Christiane Helzel, for providing Fortran code used in the shallow water example
    \item Jed Brown and Lisandro Dalcin, for help in improving serial performance
    \item Hans Petter Langtangen, for encouragement and for advice on drafts of this work
\end{itemize}
MGK acknowledges partial support from DOE Contract DE-AC01-06CH11357.
We thank the KAUST supercomputing lab for allocation of time on the
Shaheen supercomputer.  Finally, we thank the referees for many helpful
comments that led to improvements of this paper.

\bibliographystyle{siam}
\bibliography{pyclaw-sisc}

\begin{thebibliography}{10}

\bibitem{Cai2005}
{\sc X.~Cai, H.~P. Langtangen, and H.~Moe}, {\em {On the performance of the
  Python programming language for serial and parallel scientific
  computations}}, Scientific Programming, 13 (2005), pp.~31--56.

\bibitem{D.A.Calhoun2008}
{\sc D.~A. Calhoun, C.~Helzel, and R.~J. LeVeque}, {\em {Logically rectangular
  finite volume grids and methods for circular and spherical domains}}, SIAM
  Review, 50 (2008), pp.~723--752.

\bibitem{pyweno}
{\sc M.~Emmett}, {\em {PyWENO software package}}, 2011.
\newblock \url{http://memmett.github.com/PyWENO}.

\bibitem{Enkovaara2011}
{\sc J.~Enkovaara, N.~A. Romero, S.~Shende, and J.~J. Mortensen}, {\em {GPAW -
  massively parallel electronic structure calculations with Python-based
  software}}, Procedia Computer Science, 4 (2011), pp.~17--25.

\bibitem{fomel2009}
{\sc S.~Fomel and J.~F. Claerbout}, {\em {Guest Editors' Introduction:
  Reproducible Research}}, Computing in Science and Engineering, 11 (2009),
  pp.~5--7.

\bibitem{FiPy:2009}
{\sc J.~E. Guyer, D.~Wheeler, and J.~A. Warren}, {\em {FiPy: partial
  differential equations with Python}}, Computing in Science and Engineering,
  11 (2009), pp.~6--15.

\bibitem{Haurwitz1940}
{\sc B.~Haurwitz}, {\em {The motion of atmospheric disturbances on a spherical
  earth}}, Journal of Marine Research, 3 (1940), pp.~254--267.

\bibitem{ketcheson2008}
{\sc D.~I. Ketcheson}, {\em Highly efficient strong stability preserving
  {R}unge-{K}utta methods with low-storage implementations}, SIAM Journal on
  Scientific Computing, 30 (2008), pp.~2113--2136.

\bibitem{ketcheson2006}
{\sc D.~I. Ketcheson and R.~J. LeVeque}, {\em {WENOCLAW: A higher order wave
  propagation method}}, in Hyperbolic Problems: Theory, Numerics, Applications:
  Proceedings of the Eleventh International Conference on Hyperbolic Problems,
  Berlin, 2008, p.~1123.

\bibitem{sharpclaw}
{\sc D.~I. Ketcheson and M.~Parsani}, {\em {SharpClaw software}}, 2011.
\newblock \url{http://numerics.kaust.edu.sa/sharpclaw}.

\bibitem{Ketcheson2011}
{\sc D.~I. Ketcheson, M.~Parsani, and R.~J. LeVeque}, {\em {High-order wave
  propagation algorithms for general hyperbolic systems}}.
\newblock Submitted, 2011.

\bibitem{Langtangen2008}
{\sc H.~P. Langtangen and X.~Cai}, {\em {On the efficiency of Python for
  high-performance computing: A case study involving stencil updates for
  partial differential equations}}, in Modeling, Simulation and Optimization of
  Complex Processes, H.~G. Block, E.~Kostina, H.~X. Phu, and R.~Rannacher,
  eds., Springer, 2008, pp.~337--358.

\bibitem{Lax1998}
{\sc P.~Lax and X.~Liu}, {\em Solution of two-dimensional {R}iemann problems of
  gas dynamics by positive schemes}, SIAM Journal on Scientific Computing, 19
  (1998), pp.~319--340.

\bibitem{leveque1997}
{\sc R.~J. LeVeque}, {\em Wave propagation algorithms for multidimensional
  hyperbolic systems}, Journal of Computational Physics, 131 (1997),
  pp.~327--353.

\bibitem{levequefvmbook}
{\sc R.~J. LeVeque}, {\em {Finite volume methods for hyperbolic problems}},
  Cambridge University Press, Cambridge, 2002.

\bibitem{clawpack45}
{\sc R.~J. LeVeque and M.~J. Berger}, {\em {Clawpack Software version 4.5}},
  2011.
\newblock \url{http://www.clawpack.org}.

\bibitem{leveque2003}
{\sc R.~J. LeVeque and D.~H. Yong}, {\em {Solitary waves in layered nonlinear
  media}}, SIAM Journal of Applied Mathematics, 63 (2003), pp.~1539--1560.

\bibitem{pyclaw}
{\sc K.~T. Mandli, D.~I. Ketcheson, et~al.}, {\em Py{C}law software}, 2011.
\newblock Version 1.0.

\bibitem{mardal2007using}
{\sc K.~A. Mardal, O.~Skavhaug, G.~T. Lines, G.~A. Staff, and
  A.~{\O}deg{\aa}rd}, {\em Using {P}ython to solve partial differential
  equations}, Computing in Science and Engineering, 9 (2007), pp.~48--51.

\bibitem{Mortensen2005}
{\sc J.~Mortensen, L.~Hansen, and K.~Jacobsen}, {\em {Real-space grid
  implementation of the projector augmented wave method}}, Physical Review B,
  71 (2005).

\bibitem{Nilsen2010}
{\sc J.~Nilsen, X.~Cai, B.~r. H{\o}yland, and H.~P. Langtangen}, {\em
  {Simplifying the parallelization of scientific codes by a function-centric
  approach in Python}}, Computational Science \& Discovery, 3 (2010),
  p.~015003.

\bibitem{numpy}
{\sc T.~Oliphant}, {\em {NumPy Web page}}, 2011.
\newblock http://numpy.scipy.org.

\bibitem{perez2011}
{\sc F.~Perez, B.~E. Granger, and J.~D. Hunter}, {\em {P}ython: An ecosystem
  for scientific computing}, Computing in Science and Engineering, 13 (2011),
  pp.~13--21.

\bibitem{QuezadadeLuna2011}
{\sc M.~{Quezada de Luna}}, {\em Nonlinear Wave Propagation and Solitary Wave
  Formation in Two-Dimensional Heterogeneous Media}, {MSc thesis}, King
  Abdullah University of Science and Technology (KAUST), 2011.
\newblock
  \url{http://numerics.kaust.edu.sa/papers/quezada-thesis/quezada-thesis.html}.

\bibitem{quezada2011}
{\sc M.~{Quezada de Luna} and D.~I. Ketcheson}, {\em {Radial solitary waves in
  two-dimensional periodic media}}.
\newblock in preparation., 2011.

\bibitem{Shi2002}
{\sc J.~Shi, C.~Hu, and C.-W. Shu}, {\em A technique of treating negative
  weights in {WENO} schemes}, Journal of Computational Physics, 175 (2002),
  pp.~108--127.

\bibitem{Knepley2010}
{\sc V.~Stodden, M.~G. Knepley, C.~Wiggins, R.~J. LeVeque, D.~Donoho, S.~Fomel,
  M.~P. Friedlander, M.~Gerstein, I.~Mitchell, L.~L. Ouellette, N.~W. Bramble,
  P.~O. Brown, V.~Carey, L.~DeNardis, R.~Gentleman, D.~{Gezelter, J},
  A.~Goodman, J.~E. Moore, F.~A. Pasquale, J.~Rolnick, M.~Seringhaus, and
  R.~Subramanian}, {\em {Reproducible Research: addressing the need for data
  and code sharing in computational science}}, Computing in Science and
  Engineering, 12 (2010), pp.~8--13.

\bibitem{G.1968}
{\sc G.~Strang}, {\em {On the construction and comparison of difference
  schemes}}, SIAM Journal on Numerical Analysis, 5 (1968).

\bibitem{Thuburn2000}
{\sc J.~Thuburn and L.~Yong}, {\em {Numerical simulations of Rossby-Haurwitz
  waves}}, Tellus A, 52A (2000), pp.~181--189.

\bibitem{Williamson1992}
{\sc D.~L. Williamson, J.~B. Drake, J.~J. Hack, R.~Jakob, and P.~N.
  Swarztrauber}, {\em {A standard test set for numerical approximations to the
  shallow water equations on the sphere}}, Journal of Computational Physics,
  (1992), pp.~221--224.

\end{thebibliography}
\end{document}